\documentclass[12pt,a4paper,reqno]{amsart}
\usepackage{color}
\usepackage{graphics,epic}
\usepackage{amsmath,amssymb, amsthm}
\usepackage[all,2cell]{xy}

\newtheorem{theorem}{Theorem}[section]
\newtheorem*{theorem*}{Theorem}

\newtheorem{lemma}[theorem]{Lemma}
\newtheorem{proposition}[theorem]{Proposition}
\newtheorem{corollary}[theorem]{Corollary}

\newtheorem*{conjecture*}{Conjecture}

\newtheorem{remark}[theorem]{Remark}
\newtheorem{definition}[theorem]{Definition}

\renewcommand{\hat}[1]{\widehat{#1}}

\newcommand{\id}{{\rm id}}

\newcommand{\End}{{\rm End}\,}

\newcommand{\Aut}{{\rm Aut}\,}

\newcommand{\Z}{\mathbb{Z}}
\newcommand{\Q}{\mathbb{Q}}
\newcommand{\C}{\mathbb{C}}
\newcommand{\N}{\mathbb{N}}

\def\wt{{\rm wt}}
\def\C{{\mathbb C}}

\def\R{{\mathbb R}}

\def\Z{{\mathbb Z}}

\def\N{{\mathbb N}}

\def\1{{\bf 1}}
\def\tr{{\rm tr}}
\def \End{{\rm End}}

\def \Ind{{\rm Ind}}
\def \pf{\noindent {\bf Proof: \,}}
\def\theequation{5.\arabic{equation}}

\def \h{\mathfrak{h}}

\def \w{\omega}
\def \g{\mathfrak{g}}

\begin{document}

\title[Semisimplicity of module categories of certain affine VOSAs ]{Semisimplicity of module categories of certain  affine vertex operator superalgebras}
\author{Dra\v{z}en Adamovi\'{c} }
\address{Dra\v{z}en Adamovi\'{c}, Department of Mathematics, University of Zagreb, Croatia}
\thanks{ D. Adamovi\'  c was partially supported by  Croatian Science Foundation under the project IP-2022-10-9006}
\email{adamovic@math.hr}
\author{Chunrui Ai}
\address{Chunrui Ai, School of Mathematics and Statistics, Zhengzhou University, Henan 450001, China}
\email{aicr@zzu.edu.cn}

\author{Xingjun Lin}
\address{Xingjun Lin,  School of Mathematics and Statistics, Wuhan University, Wuhan 430072,  China}
\thanks{X. Lin was partially supported by China NSF grant
 12171371}
\email{linxingjun88@126.com}

\author{Jinwei Yang}
\address{Jinwei Yang, School of Mathematical Sciences, Shanghai Jiaotong University, Shanghai 200240,  China.}
\thanks{J. Yang was partially supported by China NSF grant 12371030}
\email{jinwei2@sjtu.edu.cn}

\begin{abstract}
In this paper, we show Kazhdan-Lusztig categories, that is,  the categories of lower bounded generalized weight modules for certain affine vertex operator superalgebras that  are locally finite modules of the underlying finite dimensional Lie superalgebra, are semisimple. Those are all representation categories of affine vertex operator superalgebras at conformal but non admissible levels. As a consequence, the categories of finite length generalized modules for these affine vertex operator superalgebras have braided tensor category structures.
\end{abstract}
\maketitle
\section{Introduction }
\def\theequation{1.\arabic{equation}}
\setcounter{equation}{0}
 In representation theory of vertex operator superalgebras, it is important to determine semisimplicity of certain categories of modules for a vertex operator superalgebra.  In this paper, we prove semisimplicity of Kazhdan-Lusztig categories (cf. Definition \ref{KL}) for certain simple affine vertex operator superalgebras at conformal levels that are not admissible levels. Consequently, the categories of finite length modules for these vertex operator superalgebras have braided tensor category structures.

For a finite dimensional simple Lie superalgebra $\g$ and a complex number $k$ not equal to negative dual Coxeter number of $\g$,  let $V_k(\g)$ be the simple affine vertex operator superalgebra associated to $\g$. It is well-known that $V_k(\g)$ is rational and $C_2$-cofinite if $\g$ is a Lie algebra and $k$ is a positive integer \cite{FZ}. More generally, if $k$ is an admissible level of the affine Lie algebra $\hat{\g}$ (cf. \cite{KW1}), it follows from the main result of \cite{A0} (cf.  \cite{AM},  \cite{AXL}, \cite{P1}, \cite{P2}, \cite{P3}) that Kazhdan-Lusztig category of $V_k(\g)$ is semisimple.
It remains to study Kazhdan-Lusztig category of $V_k(\g)$ if $k$ is not admissible or $\g$ is a Lie superalgebra but not a Lie algebra. The later is even more mysterious  because the category of finite dimensional modules for $\g$ is in general not semisimple.

To tackle these problems, the first named author and his collaborators used representations of minimal affine $W$-algebra $W_{k}(\g, e_{-\theta})$ to study representations of the affine vertex operator algebra $V_k(\g)$ (\cite{AKMPP3}, \cite{AP1}, \cite{AP3}, \cite{P4}). In particular, it is proved that if $k$ satisfies one of the following conditions (1)-(3),  then Kazhdan-Lusztig category of $V_k(\g)$ is semisimple.
\begin{itemize}
\item[(1)]  $k$ is a collapsing level of the minimal affine $W$--algebra $W_{k}(\g, e_{-\theta})$, that is, $W_{k}(\g, e_{-\theta})$ equals the affine vertex operator subalgebra $V_{k^{\natural}} (\g ^{\natural})$ (see Section \ref{sec:minimalw} for the details);
\item [(2)] The minimal affine $W$--algebra $W_{k}(\g, e_{-\theta})$  is a rational vertex operator algebra;
\item[(3)] The category of ordinary $W_{k}(\g, e_{-\theta})$--modules is semisimple.
\end{itemize}
This result is also extended for affine  vertex operator superalgebras \cite{AMP}. Semisimplicity of Kazhdan-Lusztig categories of some other affine vertex operator superalgebras were obtained in \cite{ACPV}, \cite{APV},  \cite{AL}, \cite{CY}. We should say that the conditions (1)-(3) are sufficient for complete reducibility. But they are not necessary conditions. In particular, the examples presented in \cite{ACPV} show that the Kazhdan-Lusztig category could be semisimple,  although the category of ordinary modules for  $W_{k}(\g, e_{-\theta})$ is not semisimple.

 If the category of ordinary modules for  $W_{k}(\g, e_{-\theta})$ is complicated, little is known for the semisimplicity of Kazhdan-Lusztig category.
In this paper, we consider affine vertex operator superalgebras appeared in conformal embeddings of minimal affine $W$-algebras \cite{AKMPP1}, \cite{AKMPP2}, and
prove Kazhdan-Lusztig categories are semisimple. Those are examples of affine vertex operator superalgebras at conformal levels. A complex number $k$ is called a conformal level if the simple minimal affine $W$-algebra $W_{k}(\g, e_{-\theta})$ and its subalgebra generated by the weight one subspace have the same conformal vectors. Conformal levels have been classified in \cite{AKMPP1}, and particularly collapsing levels are conformal levels. Concretely, we prove that Kazhdan-Lusztig categories of $V_k(\g)$ for $\g$ and $k$ listed as in Table \ref{tab1} are semisimple  (see Theorems \ref{semiconf1}, \ref{semiconf3}, \ref{semiconf2}, \ref{semiconf4}), using the results and methods of \cite{ACPV}, \cite{AKMPP2}, \cite{AKMPP3}, \cite{AMP}.

 \begin{table}[h]
\centering
\begin{tabular}{|l|l|l|}
\hline
$\g$ & $k$\\
\hline
$F(4)$ & $\frac{3}{2}$\\
\hline
$spo(2|2n), n\geq 3, (2n-1, 3)=1$& $\frac{2(n-2)}{3}$\\
\hline
$spo(2|2n+1), n\geq 1, (2n, 3)=1$ & $\frac{2n-3}{3}$\\
\hline
$sl(6|1)$ & $-2$\\
\hline
$sl(2|2n),~n\geq 2$ & $\frac{2n-1}{2}$\\
\hline
$osp(4|2n), n\geq 2$ & $n-\frac{1}{2}$\\
\hline
\end{tabular}
\caption{}
\label{tab1}
\end{table}
It is also natural to consider affine vertex operator superalgebras at collapsing levels of   arbitrary  affine $W$-algebras which are not necessarily minimal.  Explicitly, let $W_{k}(\g, f)$ be the simple affine $W$-algebra associated to $\g$ and a nilpotent element $f$ of $\g$, and let $U$ be the subalgebra of $W_{k}(\g, f)$ generated by the weight one subspace of $W_{k}(\g, f)$. A complex number $k$ is called a collapsing level of $W_{k}(\g, f)$ if $W_{k}(\g, f)=U$. Recently many collapsing levels of general affine $W$-algebras are identified in \cite{AMP}, \cite{AEM}.  However, little or nothing is known for Kazhdan-Lusztig categories of  affine vertex operator superalgebras at collapsing levels of   arbitrary  affine $W$-algebras which are not minimal.

In this paper, we consider the vertex operator algebra $V_{-2}(G_2)$ recently studied in  \cite{ADFLM} as the vertex operator algebra corresponding to rank one Argyres--Douglas theories with flavor symmetry $G_2$. Indeed, $-2$ is a collapsing level of $W_{-2}(G_2, f_{sub})$ for a subregular nilpotent element $f_{sub}$ of $G_2$ because $W_{-2}(G_2, f_{sub})=\C$ (\cite{ADFLM}). This vertex operator algebra is quasi-lisse and its irreducible ordinary modules are classified in \cite{ADFLM}. Based on this result, we prove that Kazhdan-Lusztig category of $V_{-2}(G_2)$ is semisimple (see Theorem \ref{semitypeC}). 

As a consequence of the semisimplicity results, the categories of finite length generalized modules for these affine vertex operator superalgebras have braided tensor category structures.
It is proved in \cite{CY} that the category of finite length generalized $V_k(\g)$-modules has a braided tensor category structure if all grading-restricted generalized Verma modules of $V_{k}(\g)$ are of finite length. In particular, if Kazhdan-Lusztig category of an affine vertex operator superalgebra $V_k(\g)$ is semisimple, it follows immediately that the category of finite length generalized $V_k(\g)$-modules has a braided tensor category structure.  As a consequence, we can apply abstract tensor category theory to study the category of finite length generalized $V_k(\g)$-modules. In particular,  we prove that the Kazhdan-Lusztig category of $V_{-2}(G_2)$ is equivalent to ${\rm Rep}~ S_3$ using \cite{Mc} (see Theorem \ref{rigid}). Furthermore, fusion rules of  $V_{-2}(G_2)$  are determined (see Corollary \ref{fusion}).

  We also consider the category of finite length generalized $V_k(osp(1| 2n))$-modules. It is proved recently in \cite{CGL} that the category of ordinary modules for $V_k(osp(1| 2n))$ at admissible levels is a fusion supercategory. It remains to study the category of finite length generalized modules for $V_k(osp(1| 2n))$ at non admissible levels.  In this paper,
we show that the category of finite length generalized $V_k(osp(1| 2n))$-modules has a braided tensor category structure if $k+h^\vee\notin \Q_{\geq 0}$ using \cite{CY}. This is a natural generalization of Kazhdan-Lusztig's celebrated work \cite{KL} on the existence of braided tensor categories of affine Lie algebras to the affine Lie superalgebra $\widehat{osp(1|2n)}$. In the case that $k+h^\vee\notin \Q$, it is not difficult to show that Kazhdan-Lusztig category of $V_k(osp(1| 2n))$ is semisimple. The most interesting case is when $k+h^\vee \in \Q_{<0}$.

This paper is organized as follows: In Section 2, we recall basics about vertex operator superalgebras. In Section 3, we show that Kazhdan-Lusztig categories of certain  affine vertex operator superalgebras at conformal levels are semisimple. In Section 4, we show that Kazhdan-Lusztig category of $V_{-2}(G_2)$ is semisimple. In Section 5, we show that the category of finite length generalized $V_k(osp(1| 2n))$-modules has a braided tensor category structure if $k+h^\vee\notin \Q_{\geq 0}$. We also show that Kazhdan-Lusztig category of $V_k(osp(1| 2n))$ is semisimple if $k+h^\vee\notin \Q$.

\section{Preliminaries}
\def\theequation{2.\arabic{equation}}
\setcounter{equation}{0}
\subsection{Basics on vertex operator superalgebras} In this subsection, we recall from  \cite{FHL}, \cite{KWa}, \cite{Li3} basic notions about vertex operator superalgebras.  Let $V=V_{\bar 0}\oplus V_{\bar 1}$ be a ${\Z}_2$-graded vector space, the element in $V_{\bar 0} $ (resp. $V_{\bar 1}$) is called  {\em even} (resp. {\em odd}).
We then define $[v]=i$ for any $v\in V_{\bar i}$ with  $i=0,1$.  A {\em vertex superalgebra} is a quadruple $(V,Y(\cdot, z),\1,D),$ where $V=V_{\bar 0}\oplus V_{\bar 1}$ is a ${\Z}_2$-graded vector space, ${\bf 1}$ is the {\em vacuum vector}
of $V$, $D$ is an endomorphism of $V$,   and $Y(\cdot, z)$ is a linear map
\begin{align*}
 Y(\cdot, z): V &\to (\End\,V)[[z,z^{-1}]] ,\\
 v&\mapsto Y(v,z)=\sum_{n\in{\Z}}v_nz^{-n-1}\ \ \ \  (v_n\in
\End\,V)
\end{align*}
satisfying a number of conditions (cf. \cite{Li3}).
A vertex superalgebra $V$ is called a {\em vertex operator superalgebra} if there is  a distinguished vector $\omega$, which is called the {\em conformal vector} of $V$, such that the following two conditions hold: \\
 (i) The component operators of  $Y(\w,z)=\sum_{n\in\Z}L(n)z^{-n-2}$ satisfy the Virasoro algebra
relation with {\em central charge} $c\in \C:$
\begin{align*}
[L(m),L(n)]=(m-n)L(m+n)+\frac{1}{12}(m^3-m)\delta_{m+n,0}c,
\end{align*}
and
$$L(-1)=D;$$
(ii) $V$ is $\frac{1}{2}\Z$-graded such that $V=\oplus_{n\in \frac{1}{2}\Z}V_n$,  $L(0)|_{V_n}=n$, $\dim(V_n)<\infty$ and $V_n=0$ for sufficiently small $n$. If $v\in V_n$, the {\em conformal weight} $\wt v$ of $v$ is defined to be $n$.

Next we recall from \cite{Li3}, \cite{FHL} some facts about  modules of vertex operator superalgebras.
For a  vertex operator superalgebra $V$, a {\em weak  $V$-module} $M$ is a vector space equipped
with a linear map
\begin{align*}
Y_{M}(\cdot, z):V&\to (\End M)[[z, z^{-1}]],\\
v&\mapsto Y_{M}(v,z)=\sum_{n\in\Z}v_nz^{-n-1},\,v_n\in \End M
\end{align*}
satisfying the following conditions: for any ${\Z}_2$-homogeneous  $u\in V,\ v\in V,\ w\in M$ and $n\in \Z$,
$$u_nw=0 \text{ for sufficiently large } n;$$
$$Y_M(\1, z)=\id_M;$$
\begin{align*}
\begin{split}
&z_{0}^{-1}\delta\left(\frac{z_{1}-z_{2}}{z_{0}}\right)Y_{M}(u,z_{1})Y_M(v,z_{2})-(-1)^{[u][v]}z_{0}^{-1}\delta\left(
\frac{z_{2}-z_{1}}{-z_{0}}\right)Y_M(v,z_{2})Y_M(u,z_{1})\\
&\quad=z_{2}^{-1}\delta\left(\frac{z_{1}-z_{0}}{z_{2}}\right)Y_M(Y(u,z_{0})v,z_{2}).
\end{split}
\end{align*}

A  {\em  generalized}  $V$-module  is a weak $V$-module $M$ with  a $\C$-grading  $M=\bigoplus_{\lambda\in\C}
M_{[\lambda]}$ where
$M_{[\lambda]}$ for $\lambda\in\C$ are generalized eigenspaces for the operator $L(0)$ with eigenvalues $\lambda$. A  {\em lower-bounded generalized}  $V$-module  is a generalized $V$-module $M$ such that for any $\lambda\in \C$, $M_{[\lambda+m]}=0$ for $m\in \Z$ sufficiently small.  A  {\em grading-restricted generalized}  $V$-module  is a lower-bounded generalized $V$-module $M$ such that for any $\lambda\in \C$, $\dim M_{[\lambda]}<\infty$.

For a vertex operator superalgebra $V$, an {\em ordinary} $V$-module is a grading-restricted generalized $V$-module $M$ such that $M_{[\lambda]}$ for $\lambda\in \C$ are eigenspaces for the operator $L(0)$ with eigenvalues $\lambda$.
\begin{definition}
 A generalized $V$-module $M$ is {\em finite length} if there exists generalized $V$-submodules $M=M_1\supset \cdots\supset M_{l+1}=0$ such that $M_i/M_{i+1}$ for $i=1, \cdots, l$ are irreducible ordinary $V$-modules.
\end{definition}

\subsection{Affine vertex operator superalgebras}\label{3-1}
In this subsection, we recall from \cite{FZ}, \cite{K2}, \cite{LL} some facts about  affine vertex operator superalgebras. Let $\g=\g_{\bar 0}\oplus \g_{\bar 1}$ be a basic classical simple Lie superalgebra, $\h$ be a Cartan subalgebra  of $\g_{\bar 0}$. Choose a minimal root $-\theta$ of $\g$, and let $(\cdot|\cdot)$ be the normalized  non-degenerate invariant supersymmetric bilinear form of $\g$ such that $(\theta|\theta)=2$. The affine Lie superalgebra associated to $\g$ is defined on $\hat{\g}=\g\otimes \C[t^{-1}, t]\oplus \C K$ with Lie brackets
\begin{align*}
[x(m), y(n)]&=[x, y](m+n)+(x|y) m\delta_{m+n,0}K,\\
[K, \hat\g]&=0,
\end{align*}
for $x, y\in \g$ and $m,n \in \Z$, where $x(n)$ denotes $x\otimes t^n$.

Let $M$ be a $\g$-module and $k$ be a complex number. $M$ may be viewed as a $\g\otimes \C[t]\oplus \C K$-module such that $\g\otimes t\C[t]$ acts as $0$ and $K$ acts as $k\id_{M}$. Then we consider the induced module
\begin{align*}
\hat{M}_{k}=\Ind_{\g\otimes \C[t]\oplus \C K}^{\hat \g}M.
\end{align*}

Let $ \h^*$ be the dual space of   $\h$. For $\lambda\in \h^*$, we use $L_{\g}(\lambda)$ to denote the irreducible highest weight $\g$-module of highest weight $\lambda$.
\begin{definition}\label{generV}
(1) We use $V_{\g}(k,\lambda)$ to denote the $\hat{\g}$-module $\hat{L_{\g}(\lambda)}_k$, which is called the {\em generalized Verma module}. In this paper, we also use $V^{k}(\g)$ to denote $V_{\g}(k,0)$.\\
(2) We use $L_{\g}(k,\lambda)$ to denote the irreducible quotient of  $V_{\g}(k,\lambda)$. In this paper, we also use $V_{k}(\g)$ to denote $L_{\g}(k,0)$.
\end{definition}
The following result has been established in \cite{FZ,K2}.
\begin{theorem}\label{affine}
Let $\g$ be a finite dimensional simple Lie superalgebra with the normalized invariant supersymmetric bilinear form $(\cdot|\cdot)$, $h^\vee$ be the dual Coxeter number of $(\g, (\cdot|\cdot))$ and $k$ be a complex number which is not equal to $-h^\vee$. Then $V_{k}(\g)$ is a vertex operator superalgebra with the following conformal vector:
\begin{align}\label{conformalva}
\omega=\frac{1}{2(k+h^\vee)}\sum_ia^i(-1)b^i(-1)\1,
\end{align}
where $\{a^i\}$ and $\{b^i\}$ are dual bases of $\g$ with respect to $(\cdot|\cdot)$, that is, $(b^i|a^j)=\delta_{i,j}$.
\end{theorem}

\begin{definition}\label{KL}
We will mainly focus on the following categories of $V_{k}(\g)$-modules:
\begin{itemize}
\item[(1)] Let $KL_k(\g)$ be the category of lower-bounded generalized  $V_{k}(\g)$-modules that are locally finite as $\g$-modules.
\item[(2)] Let  $KL_k^{fin}(\g)$ be the Kazhdan-Lusztig category of $V_{k}(\g)$, that is, the full subcategory of modules in $KL_k(\g)$ on which $\hat{\h}$ acts semisimply.
\item[(3)] Let $KL_k^{fl}(\g)$ be the category of finite length generalized $V_{k}(\g)$-modules.
\end{itemize}
\end{definition}

The following result has been established in Theorem 4.3 of \cite{AMP}, which  will play an important role in this paper.
\begin{theorem}\label{keysemi}
Let $\g$ be a basic Lie superalgebra. If every highest weight $V_{k}(\g)$-module in
$KL_k^{fin}(\g)$ is irreducible, then the category $KL_k^{fin}(\g)$ is semisimple.
\end{theorem}

 If $KL_k^{fin}(\g)$ is semisimple, we may show that the category $KL_k^{fl}(\g)$  has a braided tensor category structure by using the following result, which has been established in Theorem 3.10 of \cite{CY}.
\begin{theorem}\label{tensorCY}
Let $\g$ be a finite dimensional simple Lie algebra. If all grading-restricted generalized Verma modules of $V_{k}(\g)$ are of finite length, then $KL_k^{fl}(\g)$ has a braided tensor category structure.
\end{theorem}

\subsection{Heisenberg vertex operator algebras}
In this subsection, we recall some facts about
Heisenberg vertex operator algebras from \cite{FLM}, \cite{LL}. Let $H$ be a finite
dimensional vector space of dimension $d$ which has a non-degenerate symmetric bilinear
form $(\cdot|\cdot)$. Consider the Lie algebra
$$\hat{H}=H\otimes \C[t, t^{-1}]\oplus \C K$$ with Lie brackets: for $\alpha, \beta\in H$,
$$[\alpha(m), \beta(n)]=Km(\alpha| \beta)\delta_{m+n,
0},$$$$[\hat{H}, K]=0,$$
where $\alpha(n)=\alpha\otimes t^n$.

For any $\lambda\in H$, set $$M_{H}(1,
\lambda)=U(\hat{H})/J_{\lambda},$$ where $J_{\lambda}$ is the left
ideal of $U(\hat{H})$ generated by $\alpha(n),(n\geq 1)$,
$\alpha(0)-(\alpha, \lambda)$ and $K-1$. Set
$e^{\lambda}=1+J_{\lambda}$. It is  known that $M_{H}(1, \lambda)$ is
spanned by $\alpha_1(-n_1)\cdots \alpha_k(-n_k)e^{\lambda}$,
$n_1\geq \cdots \geq n_k\geq 1$. Let $\alpha_1,\cdots, \alpha_d$
be an orthonormal basis of $H$, and
 \begin{align}\label{conformalvh}
 \omega=\frac{1}{2}\sum_{1\leq
i\leq d}\alpha_i(-1)^21.
\end{align} It is known \cite{LL} that $M_{H}(1, 0)$
is a vertex operator algebra such that $\1=1+J_0$ is the vacuum vector
and $\omega$ is a conformal vector. Furthermore, $M_{H}(1,
\lambda)$ is an irreducible ordinary $M_{H}(1, 0)$-module.
 \section{Kazhdan-Lusztig categories of certain affine vertex operator superalgebras at conformal levels}
\def\theequation{3.\arabic{equation}}
\setcounter{equation}{0}
In this section, we shall prove that Kazhdan-Lusztig categories of certain  affine vertex operator superalgebras at conformal levels are semisimple.

\subsection{Minimal affine $W$-superalgebras and conformal levels}\label{sec:minimalw}
In this subsection, we shall prove that Kazhdan-Lusztig categories of certain  affine vertex operator superalgebras appeared in \cite{Lin} are semisimple. First, we recall from \cite{AKMPP2} some facts about minimal affine $W$-superalgebras and conformal levels. Let  $\g=\g_{\bar 0}\oplus \g_{\bar 1}$ be a basic classical simple Lie superalgebra, $\h$ be a Cartan subalgebra  of $\g_{\bar 0}$.  Fix a minimal root $-\theta$ of $\g$, then we choose root vectors $e_{\theta}$ and $e_{-\theta}$ such that
$$[e_{\theta}, e_{-\theta}]=x\in \h,~~[x, e_{\pm \theta}]=\pm e_{\pm \theta}.$$
Let $(\cdot|\cdot)$ be the  non-degenerate invariant supersymmetric bilinear form of $\g$ such that $(\theta|\theta)=2$, and $h^{\vee}$ be the dual Coxeter number of $(\g,(\cdot|\cdot))$. For any complex number $k\neq -h^{\vee}$, we use $W^k(\g, e_{-\theta})$ to denote the  minimal affine $W$-algebra of level $k$ defined in \cite{KRW}.  We use $\omega$ to denote the conformal vector of $W^k(\g, e_{-\theta})$ defined by the formula (2.2) of \cite{KW3}.

Let $\g^{\natural}$ be the centralizer in  $\g$ of the triple $\{e_{\theta}, x, e_{-\theta}\}$. Then it is proved in Proposition 4.1 of \cite{KRW} that $\g^{\natural}$  has  the form $$\g^{\natural}=\mathfrak{z}(\g^{\natural})\oplus \oplus_{i=1}^s\g^{\natural}_i,$$ where $\mathfrak{z}(\g^{\natural})$ denotes the center of $\g^{\natural}$ and $\g^{\natural}_1, \cdots, \g^{\natural}_s$  are simple Lie superalgebras. Furthermore, it is known  \cite{AKMPP1} that $W^k(\g, e_{-\theta})$  has a vertex subalgebra $$V^{k^{\natural}}(\g^{\natural})=M_{\mathfrak{z}(\g^{\natural})}(1, 0)\otimes \otimes_{i=1}^sV^{k^{\natural}_i}(\g^{\natural}_i)$$ for some complex numbers $k^{\natural}_1, \cdots, k^{\natural}_s$, where $M_{\mathfrak{z}(\g^{\natural})}(1, 0)$ denotes the Heisenberg vertex operator algebra associated to $\mathfrak{z}(\g^{\natural})$. Let $\omega_{\mathfrak{z}(\g^{\natural})}$, $\omega_{\g^{\natural}_i}$  be the conformal vectors of $M_{\mathfrak{z}(\g^{\natural})}(1, 0)$, $V^{k^{\natural}_i}(\g^{\natural}_i)$ defined in the formulas (\ref{conformalvh}), (\ref{conformalva}), respectively.
\begin{definition}
A complex number $k$ is called a {\em conformal level} of $\hat \g$ if $\omega=\omega_{\mathfrak{z}(\g^{\natural})}+\sum_{i=1}^s\omega_{\g^{\natural}_i}$.
\end{definition}

 When $\g$ and $k$ are as in TABLE \ref{tab2} and TABLE \ref{tab3}, it is proved in Theorem 6.8 of \cite{AKMPP2} that $k$ is a conformal level of $\hat \g$.
 \begin{table}[ht]
\centering
\begin{tabular}{|l|l|l|l|}
\hline
$\g$ & $\g^{\natural}$&  $k^{\natural}$& $k$\\
\hline
$F(4)$ &  $B_3$& $-\frac{13}{4}$ &$\frac{3}{2}$\\
\hline
$spo(2|2n), n\geq 3, (2n-1, 3)=1$ &  $D_n$& $\frac{-4n+5}{3}$& $\frac{2(n-2)}{3}$\\
 \hline
 $spo(2|2n+1), n\geq 1, (2n, 3)=1$ &  $B_n$& $\frac{-4n+3}{3}$& $\frac{2n-3}{3}$\\
 \hline
\end{tabular}
\caption{}
\label{tab2}
\end{table}

\begin{table}[ht]
\centering
\begin{tabular}{|l|l|l|l|l|l|}
\hline
$\g$ & $\g_1^\natural$& $\g_2^\natural$& $k_1^\natural$& $k_2^\natural$& $k$\\
\hline
$osp(4|2n), n\geq 2$ &  $A_1$& $C_n$& $-1/2$& $-\frac{2n+3}{4}$& $n-\frac{1}{2}$\\
 \hline
\end{tabular}
\caption{}
\label{tab3}
\end{table}

 The aim in this subsection is to show that the category $KL_k^{fin}(\g)$ is semisimple when $\g$ and $k$ are as in TABLE \ref{tab2} and TABLE \ref{tab3}. By Theorem \ref{keysemi}, it is enough to prove that every highest weight $V_{k}(\g)$-module in $KL_k^{fin}(\g)$ is irreducible. To achieve this goal, we need to recall some results from \cite{AKMPP3}. Recall that a $\widehat  \g$-module $M$ is in  category $\mathcal O ^k $    if it is $\widehat\h$-diagonalizable with finite dimensional weight spaces, $K$ acts as $k \id_M$ and  $M$ has  a finite number of maximal  weights. There is a functor  $H_\theta$ from $\mathcal O ^k $ to the category of  $W^k(\g,e_{-\theta})$-modules (see \cite{A1} for details; there $H_\theta$ is denoted by $H^0$). Let $W_k(\g, e_{-\theta})$ be the simple quotient of $W^k(\g, e_{-\theta})$. Then the following result has been established in Lemma 5.6 of \cite{AKMPP3} (see also Lemma 4.5 of \cite{AMP}).
\begin{theorem}\label{keyirr}
Let $k\in \Q\backslash \Z_{\geq 0}$. Assume that $H_\theta(U)$ is an irreducible, non-zero module of $W_k(\g,e_{-\theta})$ for every non-zero highest weight $V_k(\g)$-module $U$ from the category $KL_k^{fin}(\g)$. Then every highest weight $V_k(\g)$-module in  $KL_k^{fin}(\g)$ is irreducible.
\end{theorem}

 Therefore, to prove that every highest weight $V_{k}(\g)$-module in $KL_k^{fin}(\g)$ is irreducible, the key point is to study $H_\theta(U)$. It is proved in the proof of Theorem 4.6 of \cite{AMP} that $H_\theta(U)$ is a non-zero highest weight module for $W_k(\g,e_{-\theta})$. Moreover, we have the following result (see Section 8.2 of \cite{GK}).
 \begin{lemma}\label{ordinary}
Let $\g$ and $k$ be as in TABLE \ref{tab2} and TABLE \ref{tab3}, $U$ be a non-zero highest weight $V_{k}(\g)$-module in the category $KL_{k}^{fin}(\g)$. Then $H_\theta(U)$ is a non-zero ordinary module of $W_k(\g,e_{-\theta})$.
 \end{lemma}

 \vskip.25cm
 To show that $H_\theta(U)$ is an irreducible, non-zero $W_k(\g,e_{-\theta})$-module, we need the following results, which were established in Theorem 6.8 of \cite{AKMPP2}.
 \begin{theorem}\label{embedding1}
(1) Let $\g$, $\g^\natural$ be  simple Lie superalgebras listed in TABLE \ref{tab2}, $k, k^{\natural}$ be the complex numbers listed in TABLE \ref{tab2}. Then  $W_k(\g, e_{-\theta})$ has a subalgebra isomorphic to $V_{k^{\natural}}(\g^{\natural})$. Moreover, there exists an irreducible  $V_{k^{\natural}}(\g^{\natural})$-module $V^1$ such that $W_k(\g, e_{-\theta})=V_{k^{\natural}}(\g^{\natural})\oplus V^1$ and $\{u_{n}v|u, v\in V^1, n\in \Z\}=V_{k^{\natural}}(\g^{\natural})$.\\
(2) Let $\g$, $\g^\natural_1$,  $\g^\natural_2$ be  simple Lie superalgebras listed in TABLE \ref{tab3}, $k, k^\natural_1, k^\natural_2$ be the complex numbers listed in TABLE \ref{tab3}. Then $W_k(\g, e_{-\theta})$ has a subalgebra isomorphic to  $V_{k_1^\natural}(\g_1^\natural)\otimes V_{k_2^\natural}(\g_2^\natural)$. Moreover, there exists an irreducible  $V_{k_1^\natural}(\g_1^\natural)\otimes V_{k_2^\natural}(\g_2^\natural)$-module $V^1$ such that $W_k(\g, e_{-\theta})=V_{k_1^\natural}(\g_1^\natural)\otimes V_{k_2^\natural}(\g_2^\natural)\oplus V^1$ and $\{u_{n}v|u, v\in V^1, n\in \Z\}=V_{k_1}(\g_1^\natural)\otimes V_{k_2}(\g_2^\natural)$.
\end{theorem}

We are now ready to prove the first main result in this subsection.
\begin{theorem}\label{semiconf1}
Let $\g$ and $k$ be as in TABLE \ref{tab2}. Then the category $KL_k^{fin}(\g)$ is semisimple.
\end{theorem}
\pf By Theorem \ref{keysemi}, it is enough to prove that every highest weight $V_{k}(\g)$-module in $KL_k^{fin}(\g)$ is irreducible. Furthermore, by Theorem \ref{keyirr}, it is enough to prove that $H_\theta(U)$ is an irreducible, non-zero $W_k(\g,e_{-\theta})$-module for every non-zero highest weight $V_k(\g)$-module $U$ from the category $KL_k^{fin}(\g)$. It is proved in the proof of Theorem 4.6 of \cite{AMP} that $H_\theta(U)$ is a non-zero highest weight module for $W_k(\g,e_{-\theta})$. By Lemma \ref{ordinary}, $H_\theta(U)$ is an ordinary module of $W_k(\g,e_{-\theta})$.

We next show that $H_\theta(U)$ is an irreducible ordinary module of $W_k(\g,e_{-\theta})$. Let  $\g^\natural$ and $k^\natural$ be as in TABLE \ref{tab2}. By Theorem \ref{embedding1}, $W_k(\g, e_{-\theta})$ has a subalgebra isomorphic to $V_{k^{\natural}}(\g^{\natural})$. Moreover, there exists an irreducible  $V_{k^{\natural}}(\g^{\natural})$-module $V^1$ such that $W_k(\g, e_{-\theta})=V_{k^{\natural}}(\g^{\natural})\oplus V^1$ and $\{u_{n}v|u, v\in V^1, n\in \Z\}=V_{k^{\natural}}(\g^{\natural})$.
Let $\mathcal{C}$, $ \mathcal{D}$ be the categories of ordinary modules of $V_{k^{\natural}}(\g^{\natural})$, $W_k(\g, e_{-\theta})$, respectively. Since $k^\natural$ is an admissible level of  $\widehat{\g^\natural}$, it follows from Theorem 4.9 of \cite{CY} and Proposition 4.7 of \cite{A0} that $\mathcal{C}$ is semisimple. By Theorem \ref{embedding1}, $V_{k^{\natural}}(\g^{\natural})$ and $W_k(\g, e_{-\theta})$ have the same conformal vectors. Thus, an ordinary $W_k(\g, e_{-\theta})$-module is also an ordinary module of $V_{k^{\natural}}(\g^{\natural})$. It follows from Theorem 4.2 of \cite{Lin} that $\mathcal{D}$ is semisimple. Since $H_\theta(U)$ is an ordinary highest weight module of $W_k(\g,e_{-\theta})$, $H_\theta(U)$ must be irreducible.  This completes the proof.
\qed

\begin{remark}
If $\g$ and $k$ are as in TABLE \ref{tab2}, it is proved in Theorem 6.5 of \cite{Lin} that $W_k(\g,e_{-\theta})$ is $\omega_{\sigma}$-rational with respect to new conformal vector $\omega_{\sigma}$. However, this does not imply Theorem \ref{semiconf1} immediately. We need to prove that $H_\theta(U)$ is an admissible module for $W_k(\g,e_{-\theta})$ with respect to the new conformal vector $\omega_{\sigma}$.
\end{remark}

\vskip.25cm
We next show that $KL_k^{fin}(\g)$ is semisimple if $\g$ and $k$ are as in TABLE  \ref{tab3}. The argument is similar to that in Theorem \ref{semiconf1}. The key point is to prove that the category of ordinary $V_{-\frac{1}{2}}(A_1)\otimes V_{-\frac{2n+3}{4}}(C_n)$-modules is semisimple. In the following, we shall prove that the category of ordinary $V_{k_1}(\g_1)\otimes V_{k_2}(\g_2)$-modules is semisimple if $\g_1, \g_2$ are simple Lie algebras and $k_1, k_2$ are admissible levels of $\hat{\g}_1, \hat{\g}_2$, respectively. Note that a $V_{k_1}(\g_1)\otimes V_{k_2}(\g_2)$-module is a module of $\hat{\g}_1\oplus \hat{\g}_2$. If a $V_{k_1}(\g_1)\otimes V_{k_2}(\g_2)$-module $M$ is a highest weight module of $\hat{\g}_1\oplus \hat{\g}_2$, we say that $M$ is a highest weight module of $V_{k_1}(\g_1)\otimes V_{k_2}(\g_2)$. To prove that the category of ordinary $V_{k_1}(\g_1)\otimes V_{k_2}(\g_2)$-modules is semisimple, we need the following result.
\begin{lemma}\label{irresub1}
 Let $\g_1, \g_2$ be simple Lie algebras, $k_1, k_2$ be admissible levels of $\hat{\g}_1, \hat{\g}_2$, respectively. If $M$ is an ordinary highest weight module of  $V_{k_1}(\g_1)\otimes V_{k_2}(\g_2)$, then $M$ is an irreducible $V_{k_1}(\g_1)\otimes V_{k_2}(\g_2)$-module.
 \end{lemma}
 \pf The proof is similar to that of Proposition 2.7 of \cite{DMZ}. Let $w$ be a highest weight vector of $M$.
By the formula (4.7.4) of \cite{FHL}, $M$ is spanned by elements of the form
$$(v_{1,1}\otimes \1)_{m_1}\cdots (v_{1,s}\otimes \1)_{m_s}(\1 \otimes v_{2,1})_{n_1}\cdots (\1 \otimes v_{2,l})_{n_l}w,$$
where $v_{1,i}\in V_{k_1}(\g_1)$ and $v_{2,j}\in V_{k_2}(\g_2)$. Set $$M_1={\rm span}\{(v_{1,1}\otimes \1)_{m_1}\cdots (v_{1,s}\otimes \1)_{m_s}w|v_{1,i}\in V_{k_1}(\g_1), m_i\in \Z\}$$
and $$M_2={\rm span}\{(\1 \otimes v_{2,1})_{n_1}\cdots (\1 \otimes v_{2,l})_{n_l}w|v_{2,j}\in V_{k_2}(\g_2), n_j\in \Z\}.$$
Then  $M_1$ and $M_2$ are highest weight modules of $V_{k_1}(\g_1)$ and $V_{k_2}(\g_2)$, respectively. By the assumption, $k_1$ and $k_2$ are admissible levels of $\hat{\g}_1$ and $\hat{\g}_2$, respectively. It follows from Proposition 4.7 of \cite{A0} that $M_1$ and $M_2$ are irreducible highest weight modules of $V_{k_1}(\g_1)$ and $V_{k_2}(\g_2)$, respectively. As a consequence, $M_1\otimes M_2$ viewed as a $V_{k_1}(\g_1)\otimes V_{k_2}(\g_2)$-module is completely reducible. By the  formula (4.7.21) of \cite{FHL}, there exists a weak  $V_{k_1}(\g_1)\otimes V_{k_2}(\g_2)$-module epimorphism $\phi: M_1\otimes M_2\to M$. This implies that $M$ is an irreducible $V_{k_1}(\g_1)\otimes V_{k_2}(\g_2)$-module.
\qed

\vskip.25cm
To prove that the category of ordinary $V_{k_1}(\g_1)\otimes V_{k_2}(\g_2)$-modules is semisimple, we also need the following result.
\begin{lemma}\label{extension}
Let $L_{\g_1}(k_1, \lambda_1)\otimes L_{\g_2}(k_2, \lambda_2)$, $L_{\g_1}(k_1, \mu_1)\otimes L_{\g_2}(k_2, \mu_2)$ be irreducible ordinary highest weight modules of $V_{k_1}(\g_1)\otimes V_{k_2}(\g_2)$. If $W$ is an ordinary module of $V_{k_1}(\g_1)\otimes V_{k_2}(\g_2)$ such that
$$0\longrightarrow L_{\g_1}(k_1, \lambda_1)\otimes L_{\g_2}(k_2, \lambda_2)\longrightarrow W\longrightarrow L_{\g_1}(k_1, \mu_1)\otimes L_{\g_2}(k_2, \mu_2)\longrightarrow 0$$
is an exact sequence, then the exact sequence splits.
\end{lemma}
\pf By Lemma \ref{irresub1}, every ordinary highest weight module of  $V_{k_1}(\g_1)\otimes V_{k_2}(\g_2)$ is irreducible. If $\lambda_1\neq \mu_1$ or $\lambda_2\neq \mu_2$, we can prove that the exact sequence splits by the similar argument as that in Theorem 5.5 of \cite{AKMPP3}. If $\lambda_1=\mu_1$ and $\lambda_2=\mu_2$, let $\Omega(W)$ be the subspace of $W$ which consists of vectors having lowest conformal weight. Since $W$ is an ordinary module of  $V_{k_1}(\g_1)\otimes V_{k_2}(\g_2)$,  $\Omega(W)$ is a finite dimensional module of $\g_1\oplus \g_2$. Therefore, $\Omega(W)$ is isomorphic to $L_{\g_1}(\lambda_1)\otimes L_{\g_2}(\lambda_2)\oplus L_{\g_1}(\lambda_1)\otimes L_{\g_2}(\lambda_2)$. It follows that $W$ viewed as a $\hat{\g}_1\oplus \hat{\g}_2$-module has two highest weight vectors. By Lemma \ref{irresub1}, the exact sequence splits.
\qed

\vskip.25cm
 We now prove that the category of ordinary $V_{k_1}(\g_1)\otimes V_{k_2}(\g_2)$-modules is semisimple.
 \begin{proposition}\label{semitens}
 Let $\g_1, \g_2$ be simple Lie algebras, $k_1, k_2$ be admissible levels of $\hat{\g}_1, \hat{\g}_2$, respectively. Then the category of ordinary $V_{k_1}(\g_1)\otimes V_{k_2}(\g_2)$-modules is semisimple.
 \end{proposition}
 \pf Let $N$ be an ordinary $V_{k_1}(\g_1)\otimes V_{k_2}(\g_2)$-module, $\Omega(N)$ be the subspace of $N$ which consists of vectors having lowest conformal weight. Then $\Omega(N)$ is a finite dimensional module of $\g_1\oplus \g_2$. Therefore, $\Omega(N)$ is completely reducible. It follows that $N$ viewed as a $\hat{\g}_1\oplus \hat{\g}_2$-module contains a highest weight vector. By Lemma \ref{irresub1}, $N$ contains an irreducible highest weight module of $V_{k_1}(\g_1)\otimes V_{k_2}(\g_2)$. By Lemma \ref{extension}, the category of ordinary $V_{k_1}(\g_1)\otimes V_{k_2}(\g_2)$-modules satisfies the conditions (i)-(iv) in Lemma 1.3.1 of \cite{GK}. Therefore, the category of ordinary $V_{k_1}(\g_1)\otimes V_{k_2}(\g_2)$-modules is semisimple. This completes the proof.
\qed

\vskip.25cm
 Note that $-\frac{1}{2}$ and $-\frac{2n+3}{4}$ are admissible levels of $\hat{A_1}$ and $\hat{C_n}$, respectively. By Proposition \ref{semitens}, the category of ordinary $V_{-\frac{1}{2}}(A_1)\otimes V_{-\frac{2n+3}{4}}(C_n)$-modules is semisimple. By the similar argument as that in Theorem \ref{semiconf1}, we can prove the following result.
\begin{theorem}\label{semiconf3}
Let $\g$ and $k$ be as in TABLE  \ref{tab3}. Then the category $KL_k^{fin}(\g)$ is semisimple.
\end{theorem}

\subsection{Kazhdan-Lusztig category of  $V_{\frac{2n-1}{2}}(sl(2|2n))$} In this subsection,  we shall show that the category $KL_{\frac{2n-1}{2}}^{fin}(sl(2|2n))$ is semisimple. In the following, we will assume that  $n\geq 2$. By the proof of Theorem 8.4 of \cite{AKMPP2}, it is known that $\frac{2n-1}{2}$ is a conformal level of $\widehat{ sl(2|2n)}$. Moreover,  the vertex superalgebra $W_{\frac{2n-1}{2}}(sl(2|2n), e_{-\theta})$ has a vertex subalgebra isomorphic to $V_{-\frac{2n+1}{2}}(sl(2n))$. To study the category $KL_{\frac{2n-1}{2}}^{fin}(sl(2|2n))$, we need to recall from \cite{ACPV} some facts about modules of $V_{-\frac{2n+1}{2}}(sl(2n))$.  First, the following result has been proved in  Theorem 6.1 of \cite{ACPV}.
\begin{theorem}\label{semiACPV}
Assume that  $n\geq 2$. Then the category $KL_{-\frac{2n+1}{2}}^{fin}(sl(2n))$ is semisimple.
\end{theorem}
Let $\bar{\Lambda}_1, \cdots, \bar{\Lambda}_{2n-1}$ be the fundamental weights of $sl(2n)$.  Set
 $$U^{(2n)}_i= \left\{
\begin{array}{ll}
L_{sl(2n)}(-\frac{2n+1}{2}, i\bar{\Lambda}_1), & \text{ if } i\in \Z_{\geq 0},\\
L_{sl(2n)}(-\frac{2n+1}{2}, -i\bar{\Lambda}_{2n-1}),& \text{ if } i\in \Z_{\leq 0}.
\end{array}\right.$$
Then the following result has been established in  Theorem 6.1 of \cite{ACPV}.
\begin{theorem}\label{irrACPV}
The set $\{U^{(2n)}_i|i\in \Z\}$ is the complete list of irreducible modules in  $KL_{-\frac{2n+1}{2}}^{fin}(sl(2n))$  with fusion rules
$$U^{(2n)}_i\times U^{(2n)}_j=U^{(2n)}_{i+j}.$$
\end{theorem}

Let $\h_1=\C\alpha$ be a one dimensional vector space equipped with the bilinear form $(\alpha| \alpha)=\frac{2n}{2n-2}$. Then the following result has been established in Theorem 8.4 of \cite{AKMPP2}.
\begin{theorem}\label{conformalemb1}
 $W_{\frac{2n-1}{2}}(sl(2|2n), e_{-\theta})$ has a vertex subalgebra isomorphic to  $$V_{-\frac{2n+1}{2}}(sl(2n))\otimes M_{\h_1}(1,0).$$ Moreover, $W_{\frac{2n-1}{2}}(sl(2|2n), e_{-\theta})$ viewed as a $V_{-\frac{2n+1}{2}}(sl(2n))\otimes M_{\h_1}(1,0)$-module  has the following decomposition
 $$W_{\frac{2n-1}{2}}(sl(2|2n), e_{-\theta})=\oplus_{s\in \Z}U^{(2n)}_s\otimes M_{\h_1}\left(1,\frac{2n-2}{2n}s\alpha\right).$$
\end{theorem}

\vskip.25cm

Our aim in this subsection is to show that the category $KL_{\frac{2n-1}{2}}^{fin}(sl(2|2n))$ is semisimple. By Theorem \ref{keysemi}, it is enough to prove that every highest weight $V_{\frac{2n-1}{2}}(sl(2|2n))$-module in $KL_{\frac{2n-1}{2}}^{fin}(sl(2|2n))$ is irreducible. Furthermore, by Theorem \ref{keyirr}, it is enough to prove that $H_\theta(U)$ is an irreducible, non-zero $W_{\frac{2n-1}{2}}(sl(2|2n), e_{-\theta})$-module for every non-zero highest weight $V_{\frac{2n-1}{2}}(sl(2|2n))$-module $U$ from the category $KL_{\frac{2n-1}{2}}^{fin}(sl(2|2n))$. It is proved in the proof of Theorem 4.6 of \cite{AMP} that $H_\theta(U)$ is a non-zero highest weight module of $W_{\frac{2n-1}{2}}(sl(2|2n), e_{-\theta})$. Moreover, we have the following result (see Section 8.2 of \cite{GK}).
 \begin{lemma}\label{ordinary1}
Let $U$ be a non-zero highest weight $V_{\frac{2n-1}{2}}(sl(2|2n))$-module in the category $KL_{\frac{2n-1}{2}}^{fin}(sl(2|2n))$. Then $H_\theta(U)$ is a non-zero ordinary module of $W_{\frac{2n-1}{2}}(sl(2|2n), e_{-\theta})$.
 \end{lemma}

 \vskip.25cm
 Let $w$ be a highest weight vector of $H_\theta(U)$. Since $W_{\frac{2n-1}{2}}(sl(2|2n), e_{-\theta})$ has a vertex subalgebra isomorphic to  $V_{-\frac{2n+1}{2}}(sl(2n))\otimes M_{\h_1}(1,0),$ $H_\theta(U)$ may be viewed as a $V_{-\frac{2n+1}{2}}(sl(2n))\otimes M_{\h_1}(1,0)$-module. Let $M$ be the $V_{-\frac{2n+1}{2}}(sl(2n))\otimes M_{\h_1}(1,0)$-submodule of $H_\theta(U)$ generated by $w$. Then we have the following result.
 \begin{lemma}\label{irresub}
 Let $M$ be the $V_{-\frac{2n+1}{2}}(sl(2n))\otimes M_{\h_1}(1,0)$-submodule of $H_\theta(U)$ defined as above. Then $M$ is an irreducible $V_{-\frac{2n+1}{2}}(sl(2n))\otimes M_{\h_1}(1,0)$-module. In particular, $M$ is isomorphic to $U^{(2n)}_i\otimes M_{\h_1}(1,\mu)$ for some $i\in \Z$ and $\mu\in \h_1$.
 \end{lemma}
 \pf The proof is similar to that of Proposition 2.7 of \cite{DMZ}.  By the formula (4.7.4) of \cite{FHL}, $M$  is spanned by elements of the form
$$(v_{1,1}\otimes \1)_{m_1}\cdots (v_{1,s}\otimes \1)_{m_s}(\1 \otimes v_{2,1})_{n_1}\cdots (\1 \otimes v_{2,l})_{n_l}w,$$
where $v_{1,i}\in V_{-\frac{2n+1}{2}}(sl(2n))$ and $v_{2,j}\in M_{\h_1}(1,0)$. Set $$M_1={\rm span}\{(v_{1,1}\otimes \1)_{m_1}\cdots (v_{1,s}\otimes \1)_{m_s}w|v_{1,i}\in V_{-\frac{2n+1}{2}}(sl(2n)), m_i\in \Z\}$$
and $$M_2={\rm span}\{(\1 \otimes v_{2,1})_{n_1}\cdots (\1 \otimes v_{2,l})_{n_l}w|v_{2,j}\in M_{\h_1}(1,0), n_j\in \Z\}.$$
Then $M_1$ and $M_2$ are weak modules of $V_{-\frac{2n+1}{2}}(sl(2n))$ and $M_{\h_1}(1,0)$, respectively. Note that $M_1$ is a highest weight module of $V_{-\frac{2n+1}{2}}(sl(2n))$. Furthermore, by the similar argument as that in Proposition 4.8 of \cite{AMP}, $sl(2n)$ acts locally finitely on $M_1$. As a result, $M_1$ is an object in  $KL_{-\frac{2n+1}{2}}^{fin}(sl(2n))$. By Theorem \ref{semiACPV}, $M_1$ is an irreducible highest weight module of $V_{-\frac{2n+1}{2}}(sl(2n))$. Moreover, by Theorem \ref{irrACPV}, $M_1$ is isomorphic to $U^{(2n)}_i$ for some $i\in \Z$. Since $w$ is a  highest weight vector of $H_\theta(U)$, it follows from  Proposition 6.3.8 of \cite{LL} that $M_2$ is an irreducible module of $M_{\h_1}(1,0)$. As a consequence, $M_1\otimes M_2$ viewed as a $V_{-\frac{2n+1}{2}}(sl(2n))\otimes M_{\h_1}(1,0)$-module is completely reducible. By the  formula (4.7.21) of \cite{FHL}, there exists a weak  $V_{-\frac{2n+1}{2}}(sl(2n))\otimes M_{\h_1}(1,0)$-module epimorphism $\phi: M_1\otimes M_2\to M$. This implies that $M$ is an irreducible $V_{-\frac{2n+1}{2}}(sl(2n))\otimes M_{\h_1}(1,0)$-module. Hence, $M$ is isomorphic to $U^{(2n)}_i\otimes M_{\h_1}(1,\mu)$ for some $i\in \Z$ and $\mu\in \h_1$.
 \qed

 \vskip.25cm
 We now determine the decomposition of $H_\theta(U)$ as a $V_{-\frac{2n+1}{2}}(sl(2n))\otimes M_{\h_1}(1,0)$-module.
 \begin{lemma}\label{decomp1}
 There exist $i\in \Z$ and $\mu\in \h_1$ such that $H_\theta(U)$ has the following decomposition as a $V_{-\frac{2n+1}{2}}(sl(2n))\otimes M_{\h_1}(1,0)$-module:
 $$H_\theta(U)=\oplus_{s\in \Z}U^{(2n)}_{s+i}\otimes M_{\h_1}\left(1,\frac{2n-2}{2n}s\alpha+\mu\right).$$
 \end{lemma}
 \pf By Lemma \ref{irresub}, $M$ is isomorphic to $U^{(2n)}_i\otimes M_{\h_1}(1,\mu)$ for some $i\in \Z$ and $\mu\in \h_1$. Recall that $H_\theta(U)$ is a non-zero highest weight module of $W_{\frac{2n-1}{2}}(sl(2|2n), e_{-\theta})$ and $w$ is a  non-zero highest weight vector of $H_\theta(U)$. Since $M$ contains $w$, $H_\theta(U)$ is generated by $M$. By Proposition 4.5.6 of \cite{LL} and Theorem \ref{conformalemb1}, $H_\theta(U)$ viewed as  a $V_{-\frac{2n+1}{2}}(sl(2n))\otimes M_{\h_1}(1,0)$-module has the following decomposition
$$H_\theta(U)=\sum_{s\in \Z}U^{(2n)}_s\otimes M_{\h_1}\left(1,\frac{2n-2}{2n}s\alpha\right)\cdot M,$$
 where $U^{(2n)}_s\otimes M_{\h_1}(1,\frac{2n-2}{2n}s\alpha)\cdot M$ denotes the set $\{v_{n}c|v\in U^{(2n)}_s\otimes M_{\h_1}(1,\frac{2n-2}{2n}s\alpha), c\in M, n\in \Z\}$. By Theorem \ref{irrACPV} and the similar argument as that in Lemma 4.1 of \cite{Lin}, $U^{(2n)}_s\otimes M_{\h_1}(1,\frac{2n-2}{2n}s\alpha)\cdot M$ is an irreducible $V_{-\frac{2n+1}{2}}(sl(2n))\otimes M_{\h_1}(1,0)$-module isomorphic to  $U^{(2n)}_{s+i}\otimes M_{\h_1}(1,\frac{2n-2}{2n}s\alpha+\mu)$. Thus,
 $H_\theta(U)$ viewed as  a $V_{-\frac{2n+1}{2}}(sl(2n))\otimes M_{\h_1}(1,0)$-module has the following decomposition
$$H_\theta(U)=\oplus_{s\in \Z}U^{(2n)}_{s+i}\otimes M_{\h_1}\left(1,\frac{2n-2}{2n}s\alpha+\mu\right).$$
This completes the proof.
 \qed

 \vskip.25cm
 We are now ready to prove the main result in this subsection.
 \begin{theorem}\label{semiconf2}
 Suppose that  $n\geq 2$. Then the category $KL_{\frac{2n-1}{2}}^{fin}(sl(2|2n))$ is semisimple.
 \end{theorem}
 \pf By Theorem \ref{keysemi}, it is enough to prove that every highest weight $V_{\frac{2n-1}{2}}(sl(2|2n))$-module in $KL_{\frac{2n-1}{2}}^{fin}(sl(2|2n))$ is irreducible. Furthermore, by Theorem \ref{keyirr}, it is enough to prove that $H_\theta(U)$ is an irreducible, non-zero $W_{\frac{2n-1}{2}}(sl(2|2n), e_{-\theta})$-module for every non-zero highest weight $V_{\frac{2n-1}{2}}(sl(2|2n))$-module $U$ from the category $KL_{\frac{2n-1}{2}}^{fin}(sl(2|2n))$. It is proved in the proof of Theorem 4.6 of \cite{AMP} that $H_\theta(U)$ is a non-zero highest weight module for $W_{\frac{2n-1}{2}}(sl(2|2n), e_{-\theta})$. By Lemma \ref{ordinary1}, $H_\theta(U)$ is an ordinary module of $W_{\frac{2n-1}{2}}(sl(2|2n), e_{-\theta})$.

 We next show that $H_\theta(U)$ is an irreducible $W_{\frac{2n-1}{2}}(sl(2|2n), e_{-\theta})$-module.
Let $W$ be a non-zero $W_{\frac{2n-1}{2}}(sl(2|2n), e_{-\theta})$-submodule of $H_\theta(U)$. Then $W$  is an ordinary module of $V_{-\frac{2n+1}{2}}(sl(2n))\otimes M_{\h_1}(1,0)$. By Lemma \ref{decomp1}, $W$  viewed as  a $V_{-\frac{2n+1}{2}}(sl(2n))\otimes M_{\h_1}(1,0)$-module is completely reducible. Moreover, $W$ must contain a $V_{-\frac{2n+1}{2}}(sl(2n))\otimes M_{\h_1}(1,0)$-submodule isomorphic to $U^{(2n)}_{s+i}\otimes M_{\h_1}(1,\frac{2n-2}{2n}s\alpha+\mu)$ for some $s\in \Z$. This implies that $W$ contains the set
\begin{align}\label{subset1}
\{u_{m}b|u\in U^{(2n)}_{-s}\otimes M_{\h_1}\left(1,-\frac{2n-2}{2n}s\alpha\right), b\in U^{(2n)}_{s+i}\otimes M_{\h_1}\left(1,\frac{2n-2}{2n}s\alpha+\mu\right), m\in \Z\}.
\end{align}
By Proposition 4.5.6 of \cite{LL}, we have $\1=u_{r}v$ for some $ u\in U^{(2n)}_{-s}\otimes M_{\h_1}(1,-\frac{2n-2}{2n}s\alpha), v\in  U^{(2n)}_{s}\otimes M_{\h_1}(1,\frac{2n-2}{2n}s\alpha)$ and $r\in \Z$. Thus, for any $c\in M$, we have $$c=\1_{-1}c=(u_{r}v)_{-1}c=\sum_{i\geq 0}(-1)^i\binom{r} {i}(u_{r-i}v_{-1+i}c-(-1)^r(-1)^{[u][v]}v_{r-1-i}u_ic).$$By Proposition 4.5.8 of \cite{LL}, $v_{r-1-i}u_ic$ is a linear combination of elements of the form  $u_{p}v_{q}c$. By (\ref{subset1}), $u_{r-i}v_{-1+i}c$ and $u_{p}v_{q}c$ are elements in $W$. This implies that $c$ is an element in $W$. Therefore, $W$ contains $M$. As a consequence, $W=H_\theta(U)$. In particular, $H_\theta(U)$ is an irreducible $W_{\frac{2n-1}{2}}(sl(2|2n), e_{-\theta})$-module. This completes the proof.
 \qed
 \begin{remark}
 We also believe that the category  $KL_{k}^{fin}(sl(2|2n+1))$  is semi-simple for $k=n$, but the methods applied above can not be applied since $k \in {\Z_{\ge 1}}$ and the quantum hamiltonian reduction functor $H_{\theta}$ sends $V_k((sl(2|2n+1))$ to zero.
In the special case $n=1$,  $KL_{1}^{fin}(sl(2|3))$ is semi-simple by results presented in \cite[Section 6]{CY}.

 \end{remark}
 \subsection{Kazhdan-Lusztig category of  $V_{-2}(sl(6|1))$} In this subsection,  we shall show that the category $KL_{-2}^{fin}(sl(6|1))$  is semisimple.  By Corollary 9.3 of \cite{AKMPP2}, it is known that $-2$ is a conformal level of $\widehat{sl(6|1)}$. Moreover,  the vertex superalgebra $W_{-2}(sl(6|1), e_{-\theta})$ has a vertex subalgebra isomorphic to $V_{-1}(sl(4|1))$. To study the category $KL_{-2}^{fin}(sl(6|1))$, we need to recall from \cite{AMP}, \cite{AKMPP2} some facts about modules of $V_{-1}(sl(4|1))$.  First, the following result has been proved in  Theorem 7.7 of \cite{AMP}.
\begin{theorem}\label{semiAMP}
The category $KL_{-1}^{fin}(sl(4|1))$  is semisimple.
\end{theorem}

Let $e_{i,j}$ be the standard matrix units in $sl(4|1)$ and $$\beta=\frac{1}{3}(e_{1,1}+\cdots+e_{4,4}+4e_{5,5})_{-1}\1\in V_{-1}(sl(4|1)).$$
Following \cite{Li3}, we set
$$\Delta(\beta, z)=z^{\beta_{0}}\exp\left(\sum_{n=1}^{\infty}\frac{\beta_{n}(-z)^{-n}}{-n}\right).$$ By Proposition 5.4 of \cite{Li3},  for any $i\in \Z$,
$$\left(U_{i}, Y_{U_{i}}(\cdot, z)\right):=\left(V_{-1}(sl(4|1)), Y_{V_{-1}(sl(4|1))}(\Delta(i\beta, z)\cdot, z)\right)$$
is an irreducible weak module of $V_{-1}(sl(4|1))$. Moreover, the following result has been established in  Corollary 9.2 of \cite{AKMPP2}.
\begin{theorem}\label{irrAMKPP}
The set $\{U_i|i\in \Z\}$ is the complete list of irreducible modules in  $KL_{-1}^{fin}(sl(4|1))$  with fusion rules
$U_i\times U_j=U_{i+j}.$
\end{theorem}

Let $\h_2=\C\alpha$ be a one dimensional vector space equipped with the bilinear form $(\alpha| \alpha)=\frac{3}{5}$. Then the following result has been established in Corollary 9.3 of \cite{AKMPP2}.
\begin{theorem}\label{conformalemb2}
 $W_{-2}(sl(6|1), e_{-\theta})$ has a vertex subalgebra isomorphic to  $$V_{-1}(sl(4|1))\otimes M_{\h_2}(1,0).$$ Moreover, $W_{-2}(sl(6|1), e_{-\theta})$ viewed as a $V_{-1}(sl(4|1))\otimes M_{\h_2}(1,0)$-module  has the following decomposition
 $$W_{-2}(sl(6|1), e_{-\theta})=\oplus_{s\in \Z}U_s\otimes M_{\h_2}\left(1,\frac{5}{3}s\alpha\right).$$
\end{theorem}

\vskip.25cm

Our aim in this subsection is to show that the category $KL_{-2}^{fin}(sl(6|1))$ is semisimple. By Theorem \ref{keysemi}, it is enough to prove that every highest weight $V_{-2}(sl(6|1))$-module in $KL_{-2}^{fin}(sl(6|1))$ is irreducible. Furthermore, by Theorem \ref{keyirr}, it is enough to prove that $H_\theta(U)$ is an irreducible, non-zero $W_{-2}(sl(6|1), e_{-\theta})$-module for every non-zero highest weight $V_{-2}(sl(6|1))$-module $U$ from the category $KL_{-2}^{fin}(sl(6|1))$. It is proved in the proof of Theorem 4.6 of \cite{AMP} that $H_\theta(U)$ is a non-zero highest weight module of $W_{-2}(sl(6|1), e_{-\theta})$. Moreover, we have the following result (see Section 8.2 of \cite{GK}).
 \begin{lemma}\label{ordinary2}
Let $U$ be a non-zero highest weight $V_{-2}(sl(6|1))$-module in the category $KL_{-2}^{fin}(sl(6|1))$. Then $H_\theta(U)$ is an ordinary module of $W_{-2}(sl(6|1), e_{-\theta})$.
 \end{lemma}

 \vskip.25cm
 Let $w$ be a highest weight vector of $H_\theta(U)$. Since $W_{-2}(sl(6|1), e_{-\theta})$ has a vertex subalgebra isomorphic to  $V_{-1}(sl(4|1))\otimes M_{\h_2}(1,0)$, $H_\theta(U)$ may be viewed as a $V_{-1}(sl(4|1))\otimes M_{\h_2}(1,0)$-module. Let $M$ be the $V_{-1}(sl(4|1))\otimes M_{\h_2}(1,0)$-submodule of $H_\theta(U)$ generated by $w$. By the similar argument as that in Lemma \ref{irresub}, we have the following result.
 \begin{lemma}\label{irresub2}
 Let $M$ be the $V_{-1}(sl(4|1))\otimes M_{\h_2}(1,0)$-submodule of $H_\theta(U)$ defined as above. Then $M$ is an irreducible $V_{-1}(sl(4|1))\otimes M_{\h_2}(1,0)$-module. In particular, $M$ is isomorphic to $U_i\otimes M_{\h_2}(1,\mu)$ for some $i\in \Z$ and $\mu\in \h_2$.
 \end{lemma}

 \vskip.25cm
 By Theorems \ref{irrAMKPP}, \ref{conformalemb2} and the similar argument as that in Lemma \ref{decomp1}, we have the following decomposition of $H_\theta(U)$ as a $V_{-1}(sl(4|1))\otimes M_{\h_2}(1,0)$-module.
 \begin{lemma}\label{deomp2}
 There exist $i\in \Z$ and $\mu\in \h_2$ such that $H_\theta(U)$ has the following decomposition as a $V_{-1}(sl(4|1))\otimes M_{\h_2}(1,0)$-module:
 $$H_\theta(U)=\oplus_{s\in \Z}U_{s+i}\otimes M_{\h_2}\left(1,\frac{5}{3}s\alpha+\mu\right).$$
 \end{lemma}

 \vskip.25cm
 By Lemma \ref{deomp2} and the similar argument as that in Theorem \ref{semiconf2}, we can prove the following result.
 \begin{theorem}\label{semiconf4}
 The category $KL_{-2}^{fin}(sl(6|1))$  is semisimple.
 \end{theorem}
\section{Kazhdan-Lusztig category of $V_{-2}(G_2)$}
In this section, we shall show that the category $KL_{-2}^{fin}(G_2)$ is semisimple.  First, we recall from \cite{K} some facts about Kac-Moody algebras. Let $A=(a_{i,j})_{i,j=1}^m$ be a generalized Cartan matrix and $(\h, \Pi, \Pi^{\vee})$ be a realization of $A$ as defined in \cite{K}. In particular, $\h$ is a complex vector space, $\Pi=\{\alpha_1, \cdots, \alpha_m\}\subset \h^*$ and $\Pi^\vee=\{\alpha_1^{\vee}, \cdots, \alpha_m^{\vee}\}\subset \h$ are subsets of $\h^*$ and $\h$, respectively, such that both sets $\Pi$ and $\Pi^{\vee}$ are linearly independent.  Set  $Q_+=\{x_1\alpha_1+\cdots+x_m\alpha_m|x_i\in \Z_{\geq 0}\}$.
 Then one may define a partial order $\leq$ on $\h^*$ by $\lambda\leq \mu$ if and only if $\mu-\lambda\in Q_+$.

In the following, we let
$$A=\left(\begin{array}{lllll}
2&-1&0\\
-1&2&-1\\
0&-3&2
\end{array}\right)$$ be the generalized Cartan matrix of  type  $G_2^{(1)}$, and $\g(A)$ be the Kac-Moody algebra associated to $A$ as defined in \cite{K}. We use $e_i, f_i \ (i=0,1,2)$ to denote the Chevalley generators of $\g(A)$. Denote by $\overline{\g(A)}$ the subalgebra of $\g(A)$ generated by $e_i$ and $f_i$ with $i=1, 2$. Then it is known \cite{K} that $\overline{\g(A)}$ is a Kac-Moody algebra associated to the matrix $\bar A$ obtained from $A$ by deleting the $0$th row and column. The elements $e_i, f_i~ (i=1, 2)$ are the Chevalley generators of $\overline{\g(A)}$, and $\bar\h=\overline{\g(A)}\cap \h$ is its Cartan subalgebra. Moreover, $\overline{\g(A)}=G_2$. Let $\bar{\Lambda}_1, \bar{\Lambda}_{2}$ be the fundamental weights of $G_2$. Then the following result has been established in  Theorem 3.5 of \cite{ADFLM}.
\begin{theorem}\label{clasirr}
The set
$$\{L_{G_2}(-2, 0), L_{G_2}(-2, \bar{\Lambda}_1), L_{G_2}(-2, \bar{\Lambda}_2)\}$$
provides a complete list of irreducible ordinary $V_{-2}(G_2)$-modules.
\end{theorem}

 For any $\lambda\in \h^*$, let $M_{\g(A)}(\lambda)$ be the Verma module of $\g(A)$ of highest weight $\lambda$. We use $L_{\g(A)}(\lambda)$ to denote the simple quotient of $M_{\g(A)}(\lambda)$. Let $\delta$ be the imaginary root of $\g(A)$ defined in Theorem 5.6 of \cite{K}, $\Lambda_0\in \h^*$ be the fundamental weight of $\g(A)$ defined in Section 6.2 of \cite{K}. Then we have the following isomorphisms of $\g(A)$-modules:
\begin{align*}
&L_{G_2}(-2, 0)\cong L_{\g(A)}(-2\Lambda_0),\\
&L_{G_2}(-2, \bar{\Lambda}_1)\cong L_{\g(A)}(-2\Lambda_0+\bar{\Lambda}_1-2\delta),\\
&L_{G_2}(-2, \bar{\Lambda}_2)\cong L_{\g(A)}(-2\Lambda_0+\bar{\Lambda}_2-\delta).
\end{align*}
 In the following, we shall show that ordinary highest weight modules of $V_{-2}(G_2)$ of highest weights $-2\Lambda_0$, $-2\Lambda_0+\bar{\Lambda}_1-2\delta$ or $ -2\Lambda_0+\bar{\Lambda}_2-\delta$ are irreducible. First, we have the following result.
 \begin{proposition}\label{irreucible1}
 If $M$ is an ordinary highest weight module of $V_{-2}(G_2)$ of highest weight $-2\Lambda_0$, then $M$ is an irreducible module of $V_{-2}(G_2)$.
 \end{proposition}
 \pf Let $v$ be a highest weight vector of $M$, and $\langle v\rangle$ be the $G_2$-submodule of $M$ generated by $v$. Since  $M$ is an ordinary highest weight module of $V_{-2}(G_2)$, $\langle v\rangle$ is finite dimensional. This implies that $\langle v\rangle$ is 1-dimensional. As a result, $x(n)v$ is equal to $0$ for any $x\in G_2$ and $n\in \Z_{\geq 0}$. By the formula (3.1.12) of \cite{LL}, $v$ is a vacuum-like vector  in  $M$  (cf. Definition 4.7.1 of \cite{LL}). It follows from Proposition 4.7.7 of \cite{LL} that $M$ is isomorphic to $V_{-2}(G_2)$. Therefore, $M$ is irreducible.
 \qed

 \vskip.25cm
 Next we show that ordinary highest weight modules of $V_{-2}(G_2)$ of highest weight $-2\Lambda_0+\bar{\Lambda}_1-2\delta$ are irreducible.
 \begin{proposition}\label{irreucible2}
 If $M$ is an ordinary highest weight module of $V_{-2}(G_2)$ of highest weight $-2\Lambda_0+\bar{\Lambda}_1-2\delta$, then $M$ is an irreducible module of $V_{-2}(G_2)$.
 \end{proposition}
 \pf Assume that $N$ is a proper $V_{-2}(G_2)$-submodule of $M$. Then $N$ is a proper $\g(A)$-submodule of $M$. Viewed as a $\g(A)$-module, $N$ contains a highest weight vector $w$. Let $\lambda$ be the weight of $w$, and $\langle w\rangle$ be the $\g(A)$-submodule of $M$ generated by $w$. Then $\lambda<-2\Lambda_0+\bar{\Lambda}_1-2\delta$ and $\langle w\rangle$  is an ordinary module of $V_{-2}(G_2)$. It follows that the irreducible quotient of  $\langle w \rangle$ is an irreducible ordinary module of $V_{-2}(G_2)$. By Theorem \ref{clasirr}, $\lambda$ must be equal to $-2\Lambda_0$, $-2\Lambda_0+\bar{\Lambda}_1-2\delta$ or $ -2\Lambda_0+\bar{\Lambda}_2-\delta$. This implies that $\lambda\geq -2\Lambda_0+\bar{\Lambda}_1-2\delta$, which is a contradiction. Hence, $M$ is an irreducible module of $V_{-2}(G_2)$.
\qed

 \vskip.25cm
 In the following, for a vertex operator algebra $V$ and a generalized $V$-module $M$, we use $M_{[\lambda]}$  to denote the generalized eigenspace for the operator $L(0)$ with eigenvalue $\lambda$, where $\lambda$ is a complex number. To show that ordinary highest weight modules of $V_{-2}(G_2)$ of highest weight $ -2\Lambda_0+\bar{\Lambda}_2-\delta$ are irreducible,  we need to determine $\dim L_{G_2}(-2, \bar{\Lambda}_2)_{[2]}$. To compute $\dim L_{G_2}(-2, \bar{\Lambda}_2)_{[2]}$, we first recall from \cite{AP} some facts about the conformal embedding related to $V_{-2}(G_2)$. Let $\mathfrak{d}$ be the simple Lie algebra of type $D_4$, and $V_{-2}(\mathfrak{d})$ be the simple affine vertex operator algebra associated to $\mathfrak{d}$. Then the following result has been established in Theorem 5 of \cite{AP}.
 \begin{theorem}\label{confemb}
 $V_{-2}(\mathfrak{d})$ contains a subalgebra isomorphic to $V_{-2}(G_2)$. Moreover, we have the following isomorphism of $V_{-2}(G_2)$-modules:
 $$V_{-2}(\mathfrak{d})\cong V_{-2}(G_2)\oplus L_{G_2}(-2, \bar{\Lambda}_2)\oplus L_{G_2}(-2, \bar{\Lambda}_2)\oplus L_{G_2}(-2, \bar{\Lambda}_1).$$
 \end{theorem}
 \begin{remark}
 In this paper, $\alpha_1$ is a long root of $G_2$. Then $\alpha_1$ is equal to the second simple root of $G_2$ in \cite{ADFLM}, \cite{AP}.
 \end{remark}
 Let $c_{V_{-2}(\mathfrak{d})}$ be the central charge of $V_{-2}(\mathfrak{d})$. By direct computation, $c_{V_{-2}(\mathfrak{d})}=-14$. We use $\mathbb{H}$ to denote the complex upper half-plane, and define the normalized character $\chi_{V_{-2}(\mathfrak{d})}$ of $V_{-2}(\mathfrak{d})$ as follows:
 $$\chi_{V_{-2}(\mathfrak{d})}=\tr_{V_{-2}(\mathfrak{d})}e^{2\pi i\tau(L(0)-\frac{c_{V_{-2}(\mathfrak{d})}}{24})}\ \ \ \ \ (\tau\in \mathbb{H}).$$
 Then the following result has been proved in \cite{AK} (see Page 53 of \cite{AK}).
 \begin{theorem}\label{character}
 Let $E_{4}(\tau)$ be the normalized Eisenstein series of weight $4$, and $\eta(\tau)$ be the Dedekind eta function. Then
 $$\chi_{V_{-2}(\mathfrak{d})}=\frac{E_{4}'(\tau)}{240\eta(\tau)^{10}},$$
 where $E_{4}'(\tau)=\frac{1}{2\pi i}\frac{d}{d\tau}E_{4}(\tau)$.
 \end{theorem}

To compute $\dim L_{G_2}(-2, \bar{\Lambda}_2)_{[2]}$, we also need to determine  $\dim V_{-2}(G_2)_{[2]}$. It is proved in \cite{ADFLM} that $V^{-2}(G_2)$ has a singular vector $v_{\rm [sing]}$ of conformal weight $6$. Moreover, the following result has been proved in Theorem 3.6 of \cite{ADFLM}.
 \begin{theorem}\label{simplicity}
 Let $\langle v_{\rm [sing]}\rangle$ be the ideal of $V^{-2}(G_2)$ generated by $v_{\rm [sing]}$. Then $$V^{-2}(G_2)/\langle v_{\rm [sing]}\rangle\cong V_{-2}(G_2).$$
 \end{theorem}

As a consequence, we have the following result.
\begin{lemma}\label{dim1}
$\dim V_{-2}(G_2)_{[2]}=119$.
\end{lemma}
\pf Note that $G_2\otimes t^{-1}\C[t^{-1}]$ is a subalgebra of the affine Lie algebra $\hat{G_2}$. We use $U(G_2\otimes t^{-1}\C[t^{-1}])$ to denote the universal enveloping algebra of $G_2\otimes t^{-1}\C[t^{-1}]$. Then we have the following isomorphism of vector spaces
$$V^{-2}(G_2)\cong U(G_2\otimes t^{-1}\C[t^{-1}]).$$
By the PBW Theorem, $\dim V^{-2}(G_2)_{[2]}=119$. Since the lowest conformal weight of $\langle v_{\rm [sing]}\rangle$ is $6$, we have $\dim (V^{-2}(G_2)/\langle v_{\rm [sing]}\rangle)_{[2]}=119$. It follows from Theorem \ref{simplicity} that $\dim V_{-2}(G_2)_{[2]}=119$.
\qed

\vskip.25cm
Now we can determine the dimension of  $ L_{G_2}(-2, \bar{\Lambda}_2)_{[2]}$.
\begin{lemma}\label{dim2}
$\dim L_{G_2}(-2, \bar{\Lambda}_2)_{[2]}=98$.
\end{lemma}
\pf By Theorem \ref{character}, we have $\dim V_{-2}(\mathfrak{d})_{[2]}=329$. It follows from Theorem \ref{confemb} that
$$\dim V_{-2}(\mathfrak{d})_{[2]}=\dim V_{-2}(G_2)_{[2]}+2\dim L_{G_2}(-2, \bar{\Lambda}_2)_{[2]}+\dim L_{G_2}(-2, \bar{\Lambda}_1)_{[2]}.$$
We next compute $\dim L_{G_2}(-2, \bar{\Lambda}_1)_{[2]}.$ Recall from Definition \ref{generV} that $V_{G_2}(-2, \bar{\Lambda}_1)$ denotes the generalized Verma module $\hat{L_{G_2}(\bar{\Lambda}_1)}_{-2}$. Viewed as vector spaces, we have the following isomorphism
$$V_{G_2}(-2, \bar{\Lambda}_1)\cong U(G_2\otimes t^{-1}\C[t^{-1}])\otimes L_{G_2}(\bar{\Lambda}_1).$$
Since the lowest conformal weight of $L_{G_2}(-2, \bar{\Lambda}_1)$ is $2$, we have $$\dim L_{G_2}(-2, \bar{\Lambda}_1)_{[2]}=\dim L_{G_2}(\bar{\Lambda}_1)=14.$$
Therefore, by Lemma \ref{dim1}, $\dim L_{G_2}(-2, \bar{\Lambda}_2)_{[2]}=\frac{1}{2}(329-119-14)=98$.
\qed

\vskip.25cm
Recall from Definition \ref{generV} that $V_{G_2}(-2, \bar{\Lambda}_2)$ denotes the generalized Verma module $\hat{L_{G_2}(\bar{\Lambda}_2)}_{-2}$, we have the following result about $\dim V_{G_2}(-2, \bar{\Lambda}_2)_{[2]}$.
 \begin{lemma}\label{dim3}
$\dim V_{G_2}(-2, \bar{\Lambda}_2)_{[2]}=98$. In particular, we have
$$\dim V_{G_2}(-2, \bar{\Lambda}_2)_{[2]}=\dim L_{G_2}(-2, \bar{\Lambda}_2)_{[2]}.$$
\end{lemma}
\pf Viewed as vector spaces, we have the following isomorphism
$$V_{G_2}(-2, \bar{\Lambda}_2)\cong U(G_2\otimes t^{-1}\C[t^{-1}])\otimes L_{G_2}(\bar{\Lambda}_2).$$
Since $\dim L_{G_2}(\bar{\Lambda}_2)=7$ and $V_{G_2}(-2, \bar{\Lambda}_2)_{[1]}=L_{G_2}(\bar{\Lambda}_2)$, it follows from the PBW Theorem that $\dim V_{G_2}(-2, \bar{\Lambda}_2)_{[2]}=98$.
\qed

 \vskip.25cm
 We are now ready to show that ordinary highest weight modules of $V_{-2}(G_2)$ of highest weight $ -2\Lambda_0+\bar{\Lambda}_2-\delta$ are irreducible.
 \begin{proposition}\label{irreucible3}
 If $M$ is an ordinary highest weight module of $V_{-2}(G_2)$ of highest weight $ -2\Lambda_0+\bar{\Lambda}_2-\delta$, then $M$ is an irreducible module of $V_{-2}(G_2)$.
 \end{proposition}
 \pf Assume that $N$ is a proper $V_{-2}(G_2)$-submodule of $M$. Then $N$ is a proper $\g(A)$-submodule of $M$. Viewed as a $\g(A)$-module, $N$ contains a highest weight vector $w$. Let $\lambda$ be the weight of $w$, and $\langle w\rangle$ be the $\g(A)$-submodule of $M$ generated by $w$. Then $\langle w\rangle$  is an ordinary module of $V_{-2}(G_2)$. It follows that the irreducible quotient of  $\langle w \rangle$ is an irreducible ordinary module of $V_{-2}(G_2)$. By Theorem \ref{clasirr}, $\lambda$ must be equal to $-2\Lambda_0$, $-2\Lambda_0+\bar{\Lambda}_1-2\delta$ or $ -2\Lambda_0+\bar{\Lambda}_2-\delta$. This implies that $\lambda= -2\Lambda_0+\bar{\Lambda}_1-2\delta$. Hence, $\dim (M/N)_{[2]}<\dim M_{[2]}$.

 On the other hand, we consider the generalized Verma module $V_{G_2}(-2, \bar{\Lambda}_2)$. By Lemma \ref{dim3}, $\dim V_{G_2}(-2, \bar{\Lambda}_2)_{[2]}=98=\dim L_{G_2}(-2, \bar{\Lambda}_2)_{[2]}$. However,
$\dim V_{G_2}(-2, \bar{\Lambda}_2)_{[2]}\geq \dim M_{[2]}$ and $\dim (M/N)_{[2]}\geq \dim L_{G_2}(-2, \bar{\Lambda}_2)_{[2]}$, this implies that $\dim V_{G_2}(-2, \bar{\Lambda}_2)_{[2]}>\dim L_{G_2}(-2, \bar{\Lambda}_2)_{[2]}$, which is a contradiction. Hence, $M$ is an irreducible module of $V_{-2}(G_2)$.
\qed

\vskip.25cm
We now prove the main result in this section.
\begin{theorem}\label{semitypeC}
The category $KL_{-2}^{fin}(G_2)$ is semisimple.
\end{theorem}
\pf This follows from Propositions \ref{irreucible1}, \ref{irreucible2}, \ref{irreucible3} and Theorem \ref{keysemi}.
\qed

\vskip.25cm
By Theorem 4.9 of \cite{CY} and Theorem \ref{semitypeC}, the category $KL_{-2}^{fl}(G_2)$  is semisimple. Moreover, by Theorem \ref{tensorCY}, we have the following result.
\begin{theorem}\label{tencn}
$KL_{-2}^{fl}(G_2)$ has a braided tensor category structure.
\end{theorem}

  Finally we show that the category $KL_{-2}^{fl}(G_2)$ is rigid (cf. \cite{BK}).
 \begin{theorem}\label{rigid}
Let $S_3$ be the symmetric group of degree $3$  and ${\rm Rep}~ S_3$ be the category of finite dimensional modules of $S_3$. Then there is a braided tensor equivalence between  $KL_{-2}^{fl}(G_2)$  and ${\rm Rep}~ S_3$. In particular, the category $KL_{-2}^{fl}(G_2)$ is rigid.
\end{theorem}
\pf Let $\Aut(V_{-2}(\mathfrak{d}))$ be the automorphism group of the affine vertex operator algebra $V_{-2}(\mathfrak{d})$. By Theorem 5 of \cite{AP}, $\Aut(V_{-2}(\mathfrak{d}))$ has a subgroup isomorphic to $S_3$, and  the fixed point subalgebra $V_{-2}(\mathfrak{d})^{S_3}$ of $V_{-2}(\mathfrak{d})$  is isomorphic to $V_{-2}(G_2)$ (cf. Page 4 of \cite{ADFLM}).

It is well-known that $S_3$ has three irreducible modules $M_1, M_2, M_3$. We assume that $M_1$ is the trivial representation of $S_3$ and $M_3$ is the alternating representation of $S_3$. Then $M_2$ is the $2$-dimensional irreducible module of $S_3$. Let $\chi_i$ be the character of $M_i$, $i=1,2,3$. Then it follows from Corollary 2.5 of \cite{DLM} that $V_{-2}(\mathfrak{d})$ viewed as a $V_{-2}(\mathfrak{d})^{S_3}\otimes S_3$-module has the following decomposition
$$V_{-2}(\mathfrak{d})=V_{-2}(\mathfrak{d})^{S_3}\otimes M_1\oplus V_{-2}(\mathfrak{d})_{\chi_2}\otimes M_2\oplus V_{-2}(\mathfrak{d})_{\chi_3}\otimes M_3,$$
where $V_{-2}(\mathfrak{d})_{\chi_2}$ and $V_{-2}(\mathfrak{d})_{\chi_3}$ are irreducible modules of $V_{-2}(\mathfrak{d})^{S_3}$. Since $V_{-2}(\mathfrak{d})^{S_3}$  is isomorphic to $V_{-2}(G_2)$, it follows from Theorem \ref{confemb} that
$$V_{-2}(\mathfrak{d})_{\chi_2}\cong L_{G_2}(-2, \bar{\Lambda}_2),~~~V_{-2}(\mathfrak{d})_{\chi_3}\cong L_{G_2}(-2, \bar{\Lambda}_1).$$
Since $KL_{-2}^{fl}(G_2)$  is semisimple and $L_{G_2}(-2, 0), L_{G_2}(-2, \bar{\Lambda}_1), L_{G_2}(-2, \bar{\Lambda}_2)$ are all irreducible ordinary modules of $V_{-2}(G_2)$, then $KL_{-2}^{fl}(G_2)$ is the subcategory of the category of ordinary modules for $V_{-2}(G_2)$ generated by $L_{G_2}(-2, 0), L_{G_2}(-2, \bar{\Lambda}_1), L_{G_2}(-2, \bar{\Lambda}_2)$. By Theorem \ref{tencn}, $KL_{-2}^{fl}(G_2)$ has a braided tensor category structure. Consequently, Assumption 4.2 of \cite{Mc} holds for $V_{-2}(\mathfrak{d})^{S_3}$. By Corollary 4.8 of \cite{Mc}, there is a braided tensor equivalence between  $KL_{-2}^{fl}(G_2)$  and ${\rm Rep}~ S_3$ (see also the proof of Proposition 4.15 of \cite{Mc}). This completes the proof.
\qed

\vskip.25cm
 By Theorem \ref{rigid}, there is a braided tensor equivalence between  $KL_{-2}^{fl}(G_2)$  and ${\rm Rep}~ S_3$. As a consequence, we can determine fusion rules of $V_{-2}(G_2)$. By direct calculation on characters of $S_3$, we have the following result.
\begin{corollary}\label{fusion}
The fusion rules of $V_{-2}(G_2)$ are as follows:
\begin{align*}
& L_{G_2}(-2, \bar{\Lambda}_1)\boxtimes L_{G_2}(-2, \bar{\Lambda}_1)=L_{G_2}(-2, 0),\\
& L_{G_2}(-2, \bar{\Lambda}_2)\boxtimes L_{G_2}(-2, \bar{\Lambda}_2)=L_{G_2}(-2, 0)\oplus L_{G_2}(-2, \bar{\Lambda}_1)\oplus L_{G_2}(-2, \bar{\Lambda}_2),\\
& L_{G_2}(-2, \bar{\Lambda}_1)\boxtimes L_{G_2}(-2, \bar{\Lambda}_2)=L_{G_2}(-2, \bar{\Lambda}_2).
\end{align*}
\end{corollary}

\section{Kazhdan-Lusztig category of $V^{k}(osp(1| 2n))$}\label{sc}
\def\theequation{5.\arabic{equation}}
\setcounter{equation}{0}
In this section, we shall show that the category $KL^{fl}_k(osp(1| 2n))$ of finite length generalized $V^{k}(osp(1| 2n))$-modules has a braided tensor category structure if $k+h^\vee\notin \Q_{\geq 0}$. We shall also show that the category $KL_{k}^{fin}(osp(1| 2n))$ is semisimple if $k+h^\vee\notin \Q$.
\subsection{Basics about the Lie superalgebra $osp(1| 2n)$}\label{basic}
In this subsection,  we recall from \cite{K1}, \cite{K} some facts about the Lie superalgebra $osp(1| 2n)$.  Let $A=(a_{ij})$ be the Cartan matrix of type $B_n$. We use $I$ to denote the set $\{1, \cdots, n\}$ and $\tau$ to denote the subset  $\{n\}$ of $I$. Following \cite{K1}, we let $\tilde{G}(A, \tau)$ be the Lie superalgebra with generators $e_i, f_i, h_i, i\in I$ and the following defining relations $(i, j\in I)$:
\begin{align*}
&[e_i, f_j]=\delta_{ij}h_i,~~~~[h_i, h_j]=0,\\
&[h_i, e_j]=a_{ij}e_j, ~~~~[h_i, f_j]=-a_{ij}f_j,\\
&({\rm ad}~e_i)^{-a_{ij}+1}e_j=({\rm ad}~f_i)^{-a_{ij}+1}f_j=0,~~~\text{ if } i\neq j,\\
&[h_i]=\bar{0}, ~~~[e_i]=[f_i]=\bar{0},~~~~\text{ if } i\notin \tau,\\
&[e_i]=[f_i]=\bar{1},~~~~\text{ if } i\in \tau.
\end{align*}

 Set $\h=\C h_1\oplus \cdots \oplus \C h_n$ and define the linear functions $\alpha_i\in \h^*, i\in I,$ by $\alpha_i(h_j)=a_{ji}$. It is known \cite{K1} that $\tilde{G}(A, \tau)$ contains a unique maximal ideal $R$ for which  $R\cap \h=0$. Following \cite{K1}, we put $G(A, \tau)=\tilde{G}(A, \tau)/R$. The images in $G(A, \tau)$ of $e_i, f_i, h_i$ will be also denoted by $e_i, f_i, h_i$, respectively.
It is known \cite{K1} that $G(A, \tau)$ is isomorphic to the simple Lie superalgebra $osp(1| 2n)$.

 The Lie superalgebra $G(A, \tau)$ has a root decomposition
$$G(A, \tau)=\oplus_{\alpha\in \Delta^+}G(A, \tau)_{\alpha}\oplus \h\oplus_{\alpha\in -\Delta^+}G(A, \tau)_{\alpha},$$
where $\Delta^+$ denotes the set of positive roots of $G(A, \tau)$. We use $\Delta^+_{\bar 0}$ and $\Delta^+_{\bar 1}$ to denote the sets of positive even roots and odd roots, respectively. Define the linear function  $\rho\in \h^*$
by $\rho(h_i)=1$, $i\in I$. Then the following result has been proved in Proposition 3.1 of \cite{K1}.
\begin{proposition}\label{cweight}
$\rho=\frac{1}{2}\sum_{\alpha\in \Delta^+_{\bar 0}}\alpha-\frac{1}{2}\sum_{\alpha\in \Delta^+_{\bar 1}}\alpha$.
\end{proposition}

We now recall some facts about the bilinear form on $G(A, \tau)$. Recall \cite{K1} that a bilinear form $(\cdot|\cdot)$ on a Lie superalgebra $G=G_{\bar 0}\oplus G_{\bar 1}$ is called {\em supersymmetric} if
$$(a| b)=(-1)^{[a][b]}(b|a)$$
and {\em invariant} if
$$(a|[b,c])=([a,b]|c)~~~\text{ for any } a, b, c\in G.$$
Since $A=(a_{ij})$ is the Cartan matrix of type $B_n$, there exist an invertible diagonal matrix $D={\rm diag}~(\epsilon_1, \cdots, \epsilon_n)$ and a symmetric matrix $B=(b_{ij})$ such that $A=DB$. Moreover, $\epsilon_1, \cdots, \epsilon_n$ are positive rational numbers. Following \cite{K}, we define a symmetric bilinear form $(\cdot|\cdot)$ on $\h$ as follows:
$$(h_i|h_j)=b_{ij}\epsilon_i\epsilon_j,~~~\text{ for any } i, j\in I.$$
By the argument in the proof of Theorem 2.2 of  \cite{K},  the bilinear form $(\cdot|\cdot)$ on $\h$ can be extended to a supersymmetric invariant bilinear form on $G(A, \tau)$  (see Remark 5.4.2 of \cite{M}).
\begin{proposition}
There exists a non-degenerate supersymmetric invariant bilinear form $(\cdot|\cdot)$  on $G(A, \tau)$ such that $(h_i|h_j)=b_{ij}\epsilon_i\epsilon_j$, for any $i, j\in I$.
\end{proposition}

It is proved in Proposition 4.7 of \cite{K} that $B$ is a positive definite matrix. Set $\h_{\R}=\R h_1\oplus \cdots \oplus \R h_n$. Then we have the following result.
\begin{proposition}
Let $(\cdot|\cdot)$ be the non-degenerate supersymmetric invariant bilinear form  on $G(A, \tau)$ such that $(h_i|h_j)=b_{ij}\epsilon_i\epsilon_j$, for any $i, j\in I$. Then the restriction of $(\cdot|\cdot)$ on $\h_{\R}$ is a positive definite bilinear form.
\end{proposition}

By using $(\cdot|\cdot)$, we may define a linear isomorphism $\nu$ from $\h$ to $\h^*$ as follows:
$$\langle \nu(a), b\rangle=(a|b),~~~\text{ for any } a,b\in \h.$$
This induces a bilinear form $(\cdot|\cdot)$ on $\h^*$, which is defined by
\begin{align}\label{bilinearf1}
(\lambda|\mu)=(\nu^{-1}(\lambda)|\nu^{-1}(\mu)), ~~~\text{ for any } \lambda,\mu\in \h^*.
\end{align}
Define $\bar{\Lambda}_i, i\in I,$ by the relations
$$\bar{\Lambda}_i(h_j)=\delta_{i,j},~~~\text{ for any } i, j\in I.$$
Let $C=(c_{ij})$ be the matrix such that $c_{i,j}=(\bar{\Lambda}_i|\bar{\Lambda}_j)$. Then we have the following result.
\begin{proposition}\label{positivedef}
The matrix $C=(c_{ij})$ is positive definite.
\end{proposition}
\pf Since $(h_i|h_j)=b_{ij}\epsilon_i\epsilon_j$ and $\alpha_i(h_j)=a_{ji}=\epsilon_jb_{ji}$, we have $\nu(h_i)=\epsilon_i\alpha_i$. Since $(\cdot|\cdot)$ is a positive definite bilinear form  on $\h_{\R}$,  it follows from the formula (\ref{bilinearf1}) that  the  bilinear form $(\cdot|\cdot)$ on  $\h^*_{\R}=\R\alpha_1\oplus \cdots\oplus \R\alpha_n$ is positive definite. Note that $\alpha_i=\sum_{j\in I}\alpha_i(h_j)\bar{\Lambda}_j$, it follows that  $\bar{\Lambda}_i, i\in I$, are elements in $\h^*_{\R}$. This implies that the matrix $C=(c_{ij})$ is positive definite.
\qed

\vskip.25cm
We next recall some facts about highest weight modules of $G(A, \tau)$. In the following, we use  $\g$ to denote the Lie superalgebra $G(A, \tau)$. It is well-known \cite{K}, \cite{K1} that $\g$ has a triangular decomposition $\g=\g_+\oplus \h\oplus \g_-$, where $\g_+=\oplus_{\alpha\in \Delta^+}\g_{\alpha}$ and $\g_-=\oplus_{\alpha\in -\Delta^+}\g_{\alpha}$. For any $\Lambda\in \h^*$, we consider the induced $\g$-module
 $$V(\Lambda)=U(\g)\otimes_{U(\g_+\oplus \h)}\C,$$
 where $\C$ is viewed as a $\g_+\oplus \h$-module such that $\g_+$ acts trivially on $\C$ and $h$ acts as $\Lambda(h)\id$ for any $h\in \h$. It is well-known \cite{K}, \cite{K1} that $V(\Lambda)$ has a unique maximal proper submodule $I(\Lambda)$. We use $E^{\Lambda}$ to denote the irreducible quotient $V(\Lambda)/I(\Lambda)$. Recall \cite{K1} that  $\Lambda\in \h^*$ is called {\em dominant} if $\Lambda(h_i)\in \Z_+$ for $1\leq i\leq n-1$ and $\Lambda(h_n)\in 2\Z_+$. We use $P_+$ to denote the set of dominant weights. Then the following result has been established in Propositions 2.4, 3.4 of \cite{K1}.
 \begin{proposition}\label{finitedim}
 The simple module $E^{\Lambda}$ is finite dimensional if and only if $\Lambda\in P_+$.
 \end{proposition}

\subsection{Braided tensor category structure of $KL^{fl}_k(osp(1| 2n))$}\label{3-3}
In this subsection, we shall show that the category $KL^{fl}_k(osp(1| 2n))$ of finite length generalized $V^{k}(osp(1| 2n))$-modules has a braided tensor category structure if $k+h^\vee\notin \Q_{\geq 0}$.  We use  $\g$ and $(\cdot|\cdot)$ to denote the Lie superalgebra $osp(1| 2n)$ and the bilinear form on $osp(1| 2n)$ defined in subsection \ref{basic}, respectively. Let $h^{\vee}$ be the dual Coxeter number of $\g$ with respect to $(\cdot|\cdot)$. Then the following result has been established in Theorem 0.2.2 of \cite{GK1}.
\begin{proposition}
Suppose that $k+h^\vee\notin \Q_{\geq 0}$. Then $V^k(\g)$ is simple. In particular, $V^k(\g)=V_k(\g)$.
\end{proposition}

 To show that the category $KL^{fl}_k(\g)$ has a braided tensor category structure, we need the the following result, which can be proved by the similar argument as that in Theorem 3.10 of \cite{CY} (see also Section 2.3 of \cite{CMY}).
\begin{theorem}\label{tensorsCY}
If all grading-restricted generalized Verma modules of $V^{k}(\g)$ are of finite length, then $KL^{fl}_k(\g)$ has a braided tensor category structure.
\end{theorem}
Therefore, the key point is to prove that all grading-restricted generalized Verma modules of $V^{k}(\g)$ are of finite length. To achieve this goal,  we need the following results which follow from Theorem \ref{affine} and Propositions \ref{finitedim}, \ref{cweight} (see also \cite{LL}).
\begin{proposition}\label{cweightm}
(1) The generalized Verma module $V_{\g}(k,\lambda)$ is a weak module of $V^k(\g)$. Moreover, $V_{\g}(k,\lambda)$ is a grading-restricted generalized module of $V^k(\g)$ if and only if  $\lambda$ is dominant.\\
(2) Let $v_{\lambda}$ be a highest weight vector of $V_{\g}(k,\lambda)$, then we have  $$L(0)v_{\lambda}=c(\lambda)v_{\lambda},$$
where $c(\lambda)=\frac{(\lambda|\lambda+2\rho)}{2(k+h^{\vee})}$.
\end{proposition}

We also need the following result,  which can be proved by Proposition \ref{positivedef} and the similar argument as that in Lemma 13.2B of \cite{H}.
\begin{proposition}\label{key}
For $\Lambda\in P_+$, set $$F_{\Lambda}=\{\lambda\in P_+|(\lambda+2\rho|\lambda+2\rho)\leq (\Lambda+2\rho|\Lambda+2\rho)\}.$$
Then $F_{\Lambda}$ is a finite set.
\end{proposition}

\vskip .25cm
We are now ready to prove the main result in this subsection.
 \begin{theorem}\label{main}
Suppose that $k+h^\vee\notin \Q_{\geq 0}$. Then $KL^{fl}_k(osp(1| 2n))$ has a braided tensor category structure.
\end{theorem}
\pf By Theorem \ref{tensorsCY}, it is enough to prove that all grading-restricted generalized Verma modules of $V^{k}(\g)$ are of finite length. The argument is similar to that of Proposition 2.14 of \cite{KL}.
Let $V_{\g}(k,\Lambda)$ be a grading-restricted generalized Verma module of $V^k(\g)$. For any $V^k(\g)$-submodule $W$ of $V_{\g}(k,\Lambda)$, $W$ has the following decomposition
$$W=\oplus_{\mu\in \C} W_{[\mu]},$$
where $W_{[\mu]}$ is the generalized eigenspace of $L(0)$ with eigenvalue $\mu$. Moreover, $ W_{[\mu]}$ is finite dimensional for any $\mu\in \C$. We then define $\delta(W)=\sum_{\lambda\in F_{\Lambda}} \dim W_{[c(\lambda)]}$, where $c(\lambda)$ is defined as in Proposition \ref{cweightm}. By Proposition \ref{key}, $\delta(W)$ is well-defined.

Let  $W^1, W^2$ be $V^k(\g)$-submodules of $V_{\g}(k,\Lambda)$ such that $W^2\subsetneq W^1$. It is obvious that $\delta(W^2)\leq \delta(W^1)$. In the following, we shall prove that $\delta(W^2)<\delta(W^1)$.
 Since $W^1, W^2$ are $V^k(\g)$-submodules of $V_{\g}(k,\Lambda)$ such that $W^2\subsetneq W^1$ and $V_{\g}(k,\Lambda)$ is a grading-restricted generalized Verma module of $V^k(\g)$, $W^1/W^2$ is a non-zero grading-restricted generalized module of $V^k(\g)$. By Proposition \ref{cweightm}, there exists a dominant weight $\lambda\in P_+$ such that there is a nontrivial homomorphism $\phi: V_{\g}(k,\lambda)\to W^1/W^2$. As a result, we have $c(\lambda)\geq c(\Lambda)$ and $\dim W^1_{[c(\lambda)]}>\dim W^2_{[c(\lambda)]}$.
We next show that $\lambda\in F_{\Lambda}$. Otherwise, we have $(\lambda+2\rho|\lambda+2\rho)>(\Lambda+2\rho|\Lambda+2\rho)$, which implies that $(\lambda|\lambda+2\rho)-(\Lambda|\Lambda+2\rho)$ is a positive rational number.  On the other hand, since $c(\lambda)\geq c(\Lambda)$, it follows from Proposition \ref{cweightm} that
$$\frac{(\lambda|\lambda+2\rho)}{2(k+h^{\vee})}\geq \frac{(\Lambda|\Lambda+2\rho)}{2(k+h^{\vee})}.$$
This implies that $\frac{(\lambda|\lambda+2\rho)-(\Lambda|\Lambda+2\rho)}{2(k+h^{\vee})}\in \N$. Therefore, we have $k+h^{\vee}\in \Q_{>0}$. This is a contradiction. Then we have $\lambda\in F_{\Lambda}$. Since $\dim W^1_{[c(\lambda)]}>\dim W^2_{[c(\lambda)]}$, we have  $\delta(W^2)<\delta(W^1)$. As a consequence, all grading-restricted generalized Verma modules of $V^{k}(\g)$ are of finite length. This completes the proof. \qed

\subsection{Semisimplicity of Kazhdan-Lusztig category $KL_{k}^{fin}(osp(1| 2n))$} In this subsection, we shall show that the category $KL_{k}^{fin}(osp(1| 2n))$ is semisimple if $k+h^\vee\notin \Q$.
\begin{theorem}\label{semisC}
Suppose that  $k+h^\vee\notin \Q$. Then the category $KL_{k}^{fin}(osp(1| 2n))$ is semisimple.
\end{theorem}
\pf By Theorem \ref{keysemi}, it is enough to prove that every highest weight $V_{k}(osp(1| 2n))$-module in $KL_{k}^{fin}(osp(1| 2n))$ is irreducible. Let $M$ be a highest weight $V_{k}(osp(1| 2n))$-module in $KL_{k}^{fin}(osp(1| 2n))$ which is of highest weight $\lambda$. Assume that $M$ is not irreducible. Let $N$ be  a nontrivial $V_{k}(osp(1| 2n))$-submodule of  $M$. Then $N$ contains a highest weight $V_{k}(osp(1| 2n))$-module of highest weight $\mu$. By the proof of Theorem \ref{main}, we have
 $$\frac{(\mu|\mu+2\rho)-(\lambda|\lambda+2\rho)}{2(k+h^{\vee})}\in \N.$$This implies that $(\mu|\mu+2\rho)-(\lambda|\lambda+2\rho)=2(k+h^{\vee})s$ for some non-negative integer $s$. However, $(\mu|\mu+2\rho)-(\lambda|\lambda+2\rho)\in \Q$, this forces that $k+h^\vee\in \Q$, which is a contradiction.  Therefore, $M$ is irreducible. This completes the proof.
\qed

\vskip.5cm
{\bf Acknowledgement}. The authors wish to thank Thomas Creutzig and Robert McRae for comments on earlier drafts of this paper.

\end{document}